\documentclass[12pt]{amsart}
\usepackage{amscd,graphicx,epsfig}
\usepackage[colorlinks=true, pdfstartview=FitV, 
linkcolor=blue, citecolor=blue, urlcolor=blue]{hyperref}
\usepackage{amssymb}
\usepackage{xcolor}
\usepackage[english]{babel}
\usepackage[utf8]{inputenc}

\title[Automorphism groups and Lie algebras]
{Automorphism groups of affine varieties and their Lie algebras}
\author{Hanspeter Kraft and Mikhail Zaidenberg}

\address{Departement Mathematik und Informatik,
\newline\indent Universit\"at Basel,  
Spiegelgasse 1, CH-4051 Basel, Switzerland}
\email{hanspeter.kraft@unibas.ch}
\address{Institut Fourier, UMR 5582, 
Laboratoire de Math\'ematiques,\newline\indent
Universit\'e Grenoble Alpes, 
CS 40700, 38058 Grenoble cedex 9, France
}
\email{mikhail.zaidenberg@univ-grenoble-alpes.fr}
\date{}

\newtheorem{thm}{Theorem}[section]
\newtheorem*{thm*}{Theorem}
\newtheorem*{conj*}{Conjecture}

\newtheorem{prop}[thm]{Proposition}
\newtheorem{lem}[thm]{Lemma}
\newtheorem{cor}[thm]{Corollary}
\newtheorem*{cor*}{Corollary}

\newtheorem{que}[thm]{Question}
\newtheorem{conj}[thm]{Conjecture}

\theoremstyle{definition}
\newtheorem{defn}[thm]{Definition}
\newtheorem{exa}[thm]{Example}
\newtheorem*{exa*}{Example}

\theoremstyle{remark}
\newtheorem*{rem*}{Remark}
\newtheorem{rem}[thm]{Remark}

\theoremstyle{remark}
\newtheorem*{rems*}{Remarks}
\newtheorem{rems}[thm]{Remarks}

\newcommand{\NN}{{\mathbb N}}
\newcommand{\ZZ}{{\mathbb Z}}
\newcommand{\PP}{{\mathbb P}}
\def\TT{{\mathbb T}}
\newcommand{\FF}{{\mathbb F}}
\newcommand{\CC}{{\mathbb C}}

\newcommand{\kk}{{\mathbb K}}
\newcommand{\GG}{{\mathbb G}}

\renewcommand{\AA}{{\mathbb A}}

\newcommand{\JJJ}{\mathcal J}

\newcommand{\be}{\begin{enumerate}}
\newcommand{\ee}{\end{enumerate}}

\newcommand{\G}{\mathfrak G}
\newcommand{\V}{\mathfrak V}
\newcommand{\W}{\mathfrak W}
\renewcommand{\H}{\mathfrak H}

\newcommand{\into}{\hookrightarrow}

\newcommand{\p}{\partial}

\DeclareMathOperator{\End}{End}
\DeclareMathOperator{\id}{id}

\DeclareMathOperator{\Aut}{Aut}
\DeclareMathOperator{\SAut}{\mathfrak{U}}

\DeclareMathOperator{\Lie}{Lie}
\DeclareMathOperator{\VF}{Vec}

\DeclareMathOperator{\GL}{GL}

\DeclareMathOperator{\SL}{SL}

\DeclareMathOperator{\Ad}{Ad}

\DeclareMathOperator{\Der}{Der}
\DeclareMathOperator{\Aff}{Aff}
\DeclareMathOperator{\LND}{LND}

\DeclareMathOperator{\Spec}{Spec}

\def\Tame{\mathop{\rm Tame}}

\DeclareMathOperator{\Ga}{\GG_{\mathrm a}}
\DeclareMathOperator{\Gm}{\GG_{\mathrm m}}

\renewcommand{\subset}{\subseteq}
\renewcommand{\supset}{\supseteq}

\newcommand{\OOO}{\mathcal O}

\newcommand{\Jonq}{{\rm Jonq}}
\renewcommand{\phi}{\varphi}

\begin{document}
{\small
\begin{abstract}
This is a brief overview of a few selected chapters 
on automorphism groups of affine varieties. It
includes some open problems.
\end{abstract}
}

\maketitle
\thanks{\noindent
\textit{Mathematics Subject Classification:}
14J50, 14R20, 14L30, 14E07,\,22F50.\mbox{\hspace{11pt}}
}

{\footnotesize \tableofcontents}
\bigskip

\section{Introduction}
We provide a brief overview of 
a few selected chapters
on automorphism groups of affine varieties.
It completes the 
existing literature on the subject,
see e.g. \cite{Arz23}, \cite{Fre17}, \cite{FuKr18},  
\cite{Kal09}, \cite{KPZ17}, \cite{Kra89}, \cite{Kra96}, 
\cite{Kra17}, \cite{Lam24}, \cite{Miy24},
\cite{Sno89}, \cite{vdE00}, 
\cite{vdEKC21}, \cite{Wri17},
etc. 

Throughout the text, $X$ stands for an affine variety 
defined over an algebraically closed field $\kk$ 
of characteristic zero, 
$\Aut(X)$ stands for the group of automorphisms of $X$, 
$\AA^n$ stands for the affine $n$-space over $\kk$ and
$\Ga$ resp. $\Gm$ stands for the additive resp. 
multiplicative group 
of the field $\kk$. 

 If $\dim(X) = 1$ 
then $\Aut(X)$ is an algebraic group. 
Moreover, any finite group can be realized as 
the full automorphism group of a smooth affine curve, 
see \cite[Proposition 7.2]{Jel15}. 
The properties of $\Aut(X)$ are
ultimately related with the geometry of $X$.
For instance, if  $\Aut(X)$ is infinite, then $X$ is uniruled, i.e. 
covered by rational curves, see \cite[Theorem 1.1]{Jel15}.
Already in the case of a smooth affine surface $X$, 
the group 
$\Aut(X)$ is often infinite-dimensional and has a 
rich algebraic structure, 
see the survey article \cite{KPZ17}. 
The present survey is focused on various aspects of 
$\Aut(X)$ as an ind-group 
and on its Lie algebra. It contains a list of open problems 
partially borrowed in \cite{FuKr18} and \cite{KrZa24}.

Recall that an \emph{ind-variety} is 
an inductive limit $\V=\varinjlim A_i$ 
of algebraic varieties $A_i$ with closed embeddings 
$A_i\hookrightarrow A_{i+1}$, $i\in \NN$. 
Identifying $A_i$ with its image in $A_{i+1}$ 
we can define an ascending filtration $\V=\bigcup A_i$. 
We use both these notations alternatively 
depending on the context.
An ind-variety $\V$ is called \emph{affine}  
if all the $A_i$ are affine varieties.
It comes 
equipped with an ind-topology. Namely, $U\subset \V$ 
is open (resp. closed) if 
$U\cap A_i$ is  Zariski open 
(resp. Zariski closed) in $A_i$ for every $i=1,2,\ldots$, 
see e.g. \cite[Definition 1.1.3]{FuKr18}. 
A map $\phi\colon \V\to\V'$ of ind-varieties 
$\V=\varinjlim A_i$ 
and $\V'=\varinjlim A_i'$
is called \emph{morphism} if for any $k$ 
there is an $l$ such that 
$\phi(A_k)\subset A'_l$ 
and the induced map $A_k\to A_l'$ 
is a morphism of varieties.
A subset $U\subset\V$ 
is said to be \emph{curve-connected} 
if any pair of points $x,y\in U$
 is contained in an irreducible algebraic curve 
$C\subset A_i\cap U$ for some $i$, 
 cf. \cite[Definition 1.6.1]{FuKr18}. 
 Clearly, a curve-connected subset is connected. 

An \emph{ind-group} $\G=\varinjlim A_i$ 
is an ind-variety with 
a group structure such that the group operations 
are morphisms of ind-varieties.
An ind-group is connected if and only if it is 
curve-connected, see \cite[Remark 2.2.3]{FuKr18}. 
The connected component $\G^0$ 
of the neutral element in an ind-group $\G$ 
is an open and closed normal ind-subgroup, 
see \cite[Proposition 2.2.1]{FuKr18}. 
In the case where $\G$ is the
inductive limit of 
algebraic groups $A_i$, where 
$A_i$ is a closed subgroup of $A_{i+1}$, 
we say that $\G$ is \emph{nested}. 
This special case is considered in Section \ref{sec-4}.

The group of regular automorphisms $\Aut(X)$ 
is an affine ind-group, see e.g. 
\cite[Theorem 5.1.1]{FuKr18}. 
A closed algebraic subset of $A_i$ 
defines an algebraic family 
of automorphisms of $X$; 
such families were studied in \cite{Ram64}, 
see also  \cite[Section 3.3]{FuKr18} and \cite{Pop14}. 

A subgroup $G\subset\Aut(X)$ that is 
a closed algebraic subset of some $A_i$
has a natural structure of an algebraic group. 
It will be called an 
\emph{algebraic subgroup} of $\Aut(X)$. 
Since  $\Aut(X)$ is an affine ind-group, 
any algebraic subgroup of $\Aut(X)$ 
is an affine algebraic group.
A faithful action of a linear algebraic group $G$ on $X$ 
defines an isomorphism 
of $G$ with an algebraic subgroup of $\Aut(X)$. 

Given an ind-group $\G=\varinjlim A_i$, 
the tangent space $T_e\G=\varinjlim T_eA_i$ carries 
a natural structure of a Lie algebra denoted $\Lie(\G)$, 
see \cite[p.~189]{Sha81}, \cite[Proposition 4.2.2]{Kum02} 
and \cite[Sec. 1.9 and 2.1]{FuKr18}. 

Let $\G=\Aut(X)$ and let $\VF(X)$ 
stand for  the Lie algebra of
regular vector fields  on $X$. 
Recall that $\VF(X)$ is naturally isomorphic to 
the Lie algebra $\Der(\OOO(X))$ of derivations 
of the $\kk$-algebra $\OOO(X)$ of 
regular functions on $X$. 
In turn, $\Lie(\Aut(X))$
can be naturally identified with 
a subalgebra of $\VF(X)$, 
see \cite[Propositions 3.2.4 and 7.2.4]{FuKr18}. 
For $X=\AA^n$ this is the subalgebra 
$\VF^{\mathrm c}(\AA^n)\subset\VF(\AA^n)$ 
of all vector fields on $\AA^n$
with constant divergence, see \cite[p.~191]{Sha81} 
and \cite[Proposition 15.7.2]{FuKr18}. 
For an algebraic subgroup 
$G\subset\Aut(X)$, $\Lie(G)$ is a finite-dimensional 
Lie subalgebra of $\VF(X)$. 

Let $\G$  be an ind-group and $\H\subset\G$ be
a closed ind-subgroup. Then $\Lie(\H)$ is a Lie subalgebra 
of $\Lie(\G)$. However, it can happen 
(contrary to \cite[Theorem 1]{Sha81}) 
that $\Lie(\H)=\Lie(\G)$
while $\H$ is a proper subgroup of a  connected ind-group $\G$,
 see Example \ref{exa:FK18}. 
In this respect, the ind-groups are more complex 
than the algebraic groups.

\section{Lie algebras of vector fields}\label{sec:lavf}
\subsection{Algebraically generated subgroups 
of $\Aut(X)$ and Lie algebras} \label{ss:ags} 
Let $G\subset\Aut(X)$ be a subgroup.
One says that $G$  is 
\emph{algebraically generated} 
if it is generated by a family $\{G_i\,|\,i\in I\}$ 
of connected algebraic 
subgroups of $\Aut(X)$. 
Thus, an algebraically generated
subgroup is connected. 

Throughout this subsection 
we assume that $G$ is algebraically generated.
 We associate to $G$ 
the Lie subalgebra 
\[
L(G)=\langle \Lie (G_i)\,|
\,i\in I\rangle_{\tiny  Lie}\ \subset \VF(X)
\] 
generated by the Lie subalgebras $\Lie(G_i)$.
In fact, $L(G)$ depends only on $G$ 
and not on the choice of generating 
subgroups $G_i$, see \cite[Theorem 2.3.1(4)]{KrZa24}. 

For any ind-group $\G$
we also define the 
\emph{canonical Lie algebra} $L_{\G}$, 
see \cite[Sect.~7.1]{KrZa24}, where
\[
L_{\G} = {\rm Span}_\kk \{T_eY
\mid Y \subset \G\,\,  \text{a closed 
algebraic subset smooth in}\,\, e\}.
\]
Consider also the Zariski closure 
$\overline{G}$ in  the ind-group $\Aut(X)$. 
This is a closed ind-subgroup with Lie algebra 
$\Lie (\overline{G})$. 
Notice that $L(G)\subset L_G$ and both
are ideals in $\Lie\overline{G}$, 
see \cite[Theorem 2.3.1 and Proposition 7.1.3]{KrZa24}.
\begin{que}[{\rm \cite[Question 1]{KrZa24}}]\label{Q1}
Do we have $L(G)=L_G=\Lie (\overline{G})$?
\end{que}

The answer is affirmative 
provided $L$ is finite-dimensional.

\begin{thm}[{\rm \cite[Theorem 1]{CD03}, 
\cite[Theorem A]{KrZa24}}]\label{thmA}
Assume that $L(G)$ is finite dimensional. 
Then $G$ is an algebraic subgroup of $\Aut(X)$
and $\Lie (G) =L_G= L(G)$.
\end{thm}

\subsection{Locally finite endomorphisms 
and Lie subalgebras} 
\begin{defn} An endomorphism $\lambda$ 
of a vector space $V$ is called 
\emph{locally finite} if every $v\in V$ 
belongs to a finite-dimensional $\lambda$-invariant 
subspace of $V$. 
\end{defn} 
Every locally finite endomorphism 
$\lambda$ has a uniquely defined 
\emph{additive Jordan decomposition} 
$\lambda = \lambda_s + \lambda_n$ 
where $\lambda_s$ is semisimple, $\lambda_n$ 
locally nilpotent and 
$\lambda_s \circ \lambda_n = 
\lambda_n \circ \lambda_s$. 
Recall that a locally finite 
endomorphism $\lambda_s$ 
is \emph{semisimple} 
if the restriction of $\lambda_s$ to any 
finite-dimensional 
$\lambda_s$-invariant 
subspace of $V$  can be diagonalized. 
\begin{defn} A subspace 
$L\subset\VF(X)=\Der(\OOO(X))$ is called 
\emph{locally finite} 
if any $f\in \OOO(X)$ is contained in 
a finite dimensional $L$-invariant subspace 
of $\OOO(X)$.
\newline
Every element of a locally finite subspace 
is locally finite as an endomorphism of $\OOO(X)$ 
and thus admits a Jordan decomposition.
\end{defn}
\begin{lem}[{\rm \cite[Lemma 1.6.2]{KrZa24}}]
A locally finite subspace  $L\subset\VF(X)$ 
is finite dimensional.
\end{lem}

Let $L \subset \VF(X)$ be a Lie subalgebra 
generated by  
a family of locally finite Lie subalgebras 
$L_{i}\subset \VF(X)$, 
$i\in I$. There is the following question, 
see \cite[Question 2]{KrZa24}.
\begin{que}\label{Q2}
Is $L$ locally finite provided $L$ is finite dimensional?
\end{que}
According to Theorem \ref{thmA} the answer is positive 
if $L_i=\Lie (G_i)$, 
where the $G_i$ are algebraic subgroups of $\Aut(X)$. 
In particular, 
we have the following corollary.
\begin{cor}[{\rm \cite[Corollary of 
Theorem~E]{KrZa24}}]\label{corQ2}
Consider  a family $\{\eta_i\mid i\in I\}$
of locally nilpotent vector fields on $X$ and let
$L:=\langle \eta_i \mid i \in I\rangle_{\text{\tiny Lie}}
\subset \VF(X)$ be
the Lie algebra generated by the $\eta_i$. 
If $L$ is finite dimensional, then $L$ is locally finite.
\end{cor}
Recall that a locally nilpotent vector field 
$\eta$ on $X$ defines 
a one-dimensional algebraic subgroup 
$U=\exp(t\eta) \subset \Aut(X)$ isomorphic to $\Ga$
such that $\Lie U = \kk \,\eta$. 

A weak form of Question \ref{Q2} is as follows. 
\begin{que}\label{Q3} Let $\xi,\eta \in \VF(X)$ 
be locally finite. 
Assume that the Lie subalgebra 
$L=\langle \xi,\eta \rangle_{\text{\tiny Lie}}$ 
generated by $\xi$ and $\eta$ is finite dimensional. 
Does it follow that $L$ is locally finite? 
\end{que} 
By Corollary \ref{corQ2} the latter is true  if $\xi$ and $\eta$ 
are locally nilpotent. Indeed, under  the latter assumption
we have  by Theorem \ref{thmA} $L=\Lie(G)$ 
where $G=\langle \exp(t\xi), \exp(t\eta)\rangle$ 
is an algebraic group. 
\section{Bijective morphisms of ind-groups}
By definition, an ind-variety $\V$ comes 
equipped with a countable 
ascending filtration by algebraic varieties 
$\V=\varinjlim V_d=\bigcup_d V_d$ 
where $V_d$ is closed in $V_{d+1}$. 
Another such filtration $\V = \bigcup_e W_e$ 
is called \emph{admissible} 
if for every $n$ there is an $m$ such that 
$V_n \subset W_m$ and $W_n \subset V_m$.
There are different notions of smoothness of 
a point in an ind-variety. 
\begin{defn}[{\rm \cite[Definition 1.9.4]{FuKr18}}]
Let $\V$ be  an ind-variety, and let $x \in \V$. We say that
\begin{itemize}
\item $\V$ is \emph{strongly smooth} in $x$ 
if there is an open neighborhood of $x$ 
which has an admissible filtration consisting 
of smooth connected algebraic varieties.
\item $\V$ is \emph{geometrically smooth} in $x$ 
if there is an admissible filtration $\V=\bigcup_k V_d$ 
such that $x\in V_d$ is a smooth point for all $d$.
\end{itemize}
\end{defn}
\begin{prop} [{\rm \cite[Propositions~1.8.5, 
1.9.6~and~2.4.1]{FuKr18}}]\label{prop:smooth}
Assume that $\kk$ is uncountable. 
Let $\phi\colon \V \to \W$ 
be a bijective morphism of ind-varieties. 
\begin{itemize}
\item[(a)] 
Suppose that $\V$ is connected 
and that $\W $ has an admissible filtration 
$\W = \bigcup_k W_k$ such that all $W_k$ 
are irreducible and normal. 
Then $\phi$ is an isomorphism.
\item[(b)]
Suppose that $\V$ is curve-connected and $\W$ 
is strongly smooth in a point $w \in \W$. 
Then there is an open neighborhood $\V'$ 
of $\phi^{-1}(w)$ in $\V$ 
such that 
$\phi|_{\V'} \colon \V' \to \W$ is an open immersion. 
If $\W$ is strongly smooth in every point, then $\phi$ 
is an isomorphism.
\end{itemize}
 Let now $\phi\colon \G \to \H$ 
 be a bijective homomorphism 
 of ind-groups. 
 Assume that $\G$ is connected and $\H$ is 
 strongly smooth in $e$. 
 Then $\phi$ is an isomorphism.
\end{prop}
\begin{exa}[{\rm see \cite[Proposition~14.2.1]{FuKr18}}]
The automorphism group $\Aut(\kk\langle x,y\rangle)$ 
of the rank $2$ free associative algebra 
$\kk\langle x,y\rangle$
has a natural structure of an ind-group. 
The homomorphism of abelianization
$\kk\langle x,y\rangle\to\kk[x,y]$ induces 
a bijective (homo)morphism of ind-groups 
\[\phi\colon\Aut(\kk\langle x,y\rangle)\to 
\Aut(\kk[x,y]),\] 
see \cite{ML70}, \cite{Cze71-72} and 
\cite[Theorem~9.3]{Coh85}.
The induced surjective morphism of Lie algebras 
\[\Lie(\Aut(\kk\langle x,y\rangle))\to \Lie(\Aut(\kk[x,y]))\]
has a non-trivial kernel, see 
\cite[last paragraph of sect. 11]{BW00} 
and \cite[Proposition 14.2.1]{FuKr18}. In particular, 
$\phi$ is not an isomorphism
of ind-groups.
In fact, there exists a unicuspidal curve $C$ in 
$\Aut(\kk\langle x,y\rangle)$ such that 
$\phi(C)\subset \Aut(\kk[x,y])$ 
has a deeper cusp and so, is not isomorphic 
to $C$, see \cite[Sect. 14.3]{FuKr18}.
\end{exa}
Using Proposition \ref{prop:smooth} 
one arrives at the following conclusion.
\begin{cor}[{\rm \cite[Corollary 14.1.2]{FuKr18}}] 
Let $\kk$ be uncountable. 
Then the ind-group $\G= \Aut(\kk[x,y])$ 
is not strongly smooth in $\id$. 
Moreover, there is
no admissible filtration 
$\G= \cup_d \G_d$ with irreducible and 
normal algebraic varieties $\G_d$. 
\end{cor}
\begin{rem} It is not clear how the notions 
``strongly smooth'' and ``geometrically smooth'' are related. 
There exists as well a notion of \emph{algebraic smoothness} 
of a point in an ind-variety.
It has the following advantage: the point $e\in \G$ is 
algebraically smooth in every affine ind-group $\G$,
see \cite[Theorem 3]{Sha81} and \cite[Theorem 4.3.7]{Kum02}. 
A geometrically smooth point of an ind-variety is 
algebraically smooth, see \cite[Lemma 4.3.4]{Kum02}.

There are several pathological examples. 
For instance, there is an ind-variety $\V$ 
such that each point $v\in \V$ is not algebraically 
smooth, see \cite[P. 188]{Sha81}. The ind-group 
 $\G=\SL(2, \kk[t])$ has an admissible 
ascending filtration $\G=\cup_d G_d$ 
such that the neutral element $e$ 
of $\G$ is singular in every $G_d$, 
see \cite[Theorem 2]{Sha04} 
and \cite[Examples 4.3.5 and 4.3.8]{Kum02}. 
\end{rem}
\section{Nested ind-groups}\label{sec-4}
\begin{defn} A subgroup $G$ of an ind-group $\G$ is 
\emph{nested} if $G$ is 
a union of an increasing sequence of algebraic 
subgroups $G_i$ of $\G$, 
see \cite{KPZ17}. 
A nested ind-group $\G$ is called \emph{unipotent} if all 
the $G_i$ are unipotent algebraic groups. 
\end{defn}
The last assertion of Proposition \ref{prop:smooth} 
can be applied to nested ind-groups. 
Indeed, such a group is strongly smooth, 
see \cite[Example 2.4.2]{FuKr18}. 
Notice that a connected nested ind-group 
is algebraically generated, 
and any closed ind-subgroup of a nested 
ind-group is nested.

\begin{que}[\rm{\cite[Question 4]{KrZa24}}] 
Is it true that every connected nested 
ind-subgroup $G\subset\Aut(X)$ is closed in $\Aut(X)$?
\end{que}
The answer is affirmative, see 
\cite[Theorem 6.6]{Per24}.
Furthermore, a subgroup of $\Aut(X)$ 
consisting of 
unipotent elements is nested 
if and only if it is closed in $\Aut(X)$,
see  \cite[Theorem 8.11]{Per24}. 
The proofs of these results exploit 
a new notion of a de Jonqui\`eres-like 
(dJ-like, for short) subgroup of $\Aut(X)$.
The following analogue of the Lie--Kolchin 
theorem holds: 
every nested unipotent subgroup 
$U \subset \Aut(X)$ 
is contained in a dJ-like subgroup, 
see  \cite[Theorem 6.1]{Per24}.
Yet another ingredient in the proof
 of \cite[Theorem 6.6]{Per24} is
the following version of the Levi decomposition 
for nested subgroups, see  
\cite[Theorem 2.11 and Corollary 2.12]{KPZ17}.
\begin{thm} Let $G=\varinjlim_i G_i$ 
be a connected nested 
subgroup of $\Aut(X)$. Then there is a decomposition 
$G = R_{\mathrm u}(G)\rtimes L$, 
where $L$ is a maximal reductive algebraic 
subgroup in $G$ 
and $R_{\mathrm u}(G)$ is 
the unipotent radical of $G$. 
Moreover, one may assume that 
$G_i = R_{\mathrm u}(G_i)\rtimes L$ 
where $R_{\mathrm u}(G_i)= R_{\mathrm u}(G)\cap G_i$.
\end{thm} 
{There is the following conjectural  criterion of nestedness. 
\begin{conj}[{\rm cf. \cite[9.4.3-9.4.6]{FuKr18}}] 
Assume that every finite subset of an ind-group 
$\G$ is contained in 
an algebraic group. Then $\G$ 
  is nested.
\end{conj}
Recall the following definition.  
\begin{defn}[{\rm\cite[Definition 1.13.2]{FuKr18}}] 
A subset $S$ of an ind-group $\G$ is called 
\emph{weakly closed} 
if for every algebraic subset
$U \subset\G$ such that $U \subset S$ 
we have $\overline{U} \subset S$. 
The union $\cup_U \overline{U}$ for all $U$ 
as before is called the \emph{weak closure} of $S$. 
\end{defn}
Notice that a countable union of 
closed algebraic subsets in $\G$ 
is weakly closed, see \cite[Proposition 1.13.6]{FuKr18}.
\begin{que}[{\rm \cite[Question 15.9.1]{FuKr18}}]
\label{que:conj}
Let $g \in \Aut(\AA^n)$. Are the
following assertions equivalent?
\begin{itemize}
\item $g$ is diagonalizable.
\item $g$ is semisimple.
\item The conjugacy class $C(g)$ 
is closed in $\Aut(\AA^n)$.
\item The conjugacy class $C(g)$ 
is weakly closed in $\Aut(\AA^n)$.
\end{itemize}
\end{que}
The next proposition contains a partial 
answer to Question \ref{que:conj}.
 \begin{prop}[{\rm \cite[Corollary 15.9.8]{FuKr18}}]
Let $g \in \Aut(\AA^n)$ be semisimple. 
Then $g$ is diagonalizable if and only if its 
conjugacy class $C(g)$ is weakly closed.
\end{prop}
\begin{que}[{\rm \cite[Question 15.9.11]{FuKr18};  
cf. \cite[Proposition 15.9.12]{FuKr18}}]
Let $g \in\Aut(\AA^n)$ be locally finite. 
Does the (weak) closure of $C(g)$ 
contain the semisimple part of $g$?
\end{que}
\section{Algebraic subgroups and algebraic elements}
\begin{defn} An element $g$ of an ind-group 
$\G$ is called \emph{algebraic} 
if it is contained in an algebraic subgroup of $\G$. 
It is called \emph{unipotent} if it is contained in 
a $\Ga$-subgroup of $\G$. 
\end{defn}
Clearly, $g$ is algebraic iff the closure 
$\overline{<g>}$ is an algebraic group.
If $g$ is algebraic, then this abelian group 
has the form $(\Gm)^s (\Ga)^t F$  
with $s \ge 0$, $t \in \{0,1\}$ and $F$ 
being a finite cyclic group. 
Every subgroup of this form in $\G$ is equal to 
$\overline{<g>}$ for some $g\in\G$.

It is easily seen that $g\in\Aut(X)$ 
is algebraic if and only if 
it defines a locally finite automorphism of $\OOO(X)$. 
In particular, any algebraic element $g$ has 
a unique \emph{multiplicative Jordan decomposition} 
$g=g_{\mathrm s}g_{\mathrm u}$, 
where $g_{\mathrm s}, g_{\mathrm u}\in \Aut(X)$ commute,  
$g_{\mathrm s}$ is semisimple and 
$g_{\mathrm u}$ is unipotent, see \cite[Section 9.1]{FuKr18}.

Let $\Aut(X)_{\mathrm {alg}}\subset\Aut(X)$ stand 
for the subgroup generated by all algebraic elements,
in other words, by all algebraic subgroups of $\Aut(X)$.
\begin{thm}[{\rm\cite[Theorem 1.1]{BD15}}] 
Assume that the field $\kk$ is uncountable. 
Let $S$ be an affine surface given in 
$\AA^4=\Spec(\kk[x,y,z,u])$
by the equations
\[yu=xP(x),\quad vx=uQ(u),\quad
yv = P(x)Q(u),\]
where each of the polynomials $P,Q\in\kk[w]$ 
has at least 2 distinct roots and $P (0)\neq 0$. 
Then the following hold.
\begin{itemize}
\item The normal subgroup 
$\Aut(S)_{\mathrm {alg}} \subset \Aut(S)$ 
is not generated by a countable family 
of algebraic subgroups, and
\item the quotient $\Aut(S)/\Aut(S)_{\mathrm {alg}}$ 
contains a free group over an uncountable set
of generators. 
\end{itemize}
\end{thm}
\begin{conj}
Let $G =\langle G_1,\ldots,G_n\rangle$ 
be a subgroup of $\Aut(X)$ 
generated by a finite collection of algebraic 
subgroups $G_i$. 
Then either $G$ is an algebraic group, or $G$ 
contains non-algebraic elements.
\end{conj}
A weak version:
\begin{conj} Let $g,h \in \Aut(X)$ be algebraic elements. 
Then either $g,h$ are contained in an 
algebraic subgroup of $\Aut(X)$, 
or the group generated by $g$ and $h$ 
contains non-algebraic elements.
\end{conj}
The following corollary of \cite[Theorem 5.1]{LPU23}
confirms the latter conjecture in dimension two.
\begin{thm}
Let $X$ be an affine surface 
and $G$ be a finitely generated subgroup of $\Aut(X)$. If $G$
consists of algebraic
elements,
then $G$ is contained in an algebraic subgroup of $\Aut(X)$.
\end{thm} 
Cf. also Theorem \ref{thm:PR22} and Corollary \ref{cor:PR24}
below.

Due to the next results (see \cite[Theorem 1.1  
and Corollaries 1.2-1.3]{BPZ24}),  the Burnside problem 
has a positive solution in $\Aut(X)$.
\begin{thm} Let $k$ be any field of zero characteristic 
and $X$ be an affine variety over $k$. If
$G$ is a torsion subgroup of $\Aut(X)$, 
then any finitely generated subgroup of $G$ is finite.
If $X=\AA^n_k$, then $G$ contains an abelian 
normal subgroup of finite index. 
\end{thm}
The proof uses the following result of 
Bass and Lubotzky \cite[Corollary 1.2]{BL83}.
\begin{thm} Let $k$ be an arbitrary field and $X$ 
be a scheme of finite type over $k$.
If $G\subset\Aut(X)$ is a finitely generated subgroup,
then $G$ is residually finite, i.e. the intersection of 
all finite index subgroups of $G$ is trivial. If $X$ is 
flat over $\ZZ$, then $G$ is virtually torsion free,
in particular, the orders of finite subgroups of 
$G$ are bounded above.
\end{thm}
\begin{exa}[\rm{see \cite[Theorem 6.4.2]{KrZa24}}]
Let $F=\langle g,h\rangle \subset\Aut(\AA^2)$
be the subgroup generated by algebraic elements
\[g\colon (x,y)\mapsto (x+y^2,y)\quad\text{and}\quad 
h\colon (x,y)\mapsto (x,y+x^2).\]
Then $F\cong\FF_2$ is a free group of rank 2. 
It contains non-algebraic elements, for example, 
$gh$ is one of them. 
The closure  $\mathcal{F}=\overline{F}$ is 
a free product $\mathcal{F}=\JJJ^+*\JJJ^-$ 
of two nested 
unipotent closed subgroups $\JJJ^\pm$ of 
$\Aut(\AA^2)$ isomorphic to 
the vector group $(\kk[t],+)$. Every algebraic 
subgroup of $\mathcal{F}$ 
is conjugate to a subgroup of $\JJJ^+$ or $\JJJ^-$.
\end{exa}
\begin{conj}[{\cite[Question~9.1.5]{FuKr18}}]\label{conj-6}
Assume that $\Aut^0(X)\neq\{\id_X\}$.
Then $\Aut^0(X)$ contains 
a subgroup 
isomorphic to either $\mathbb{G}_{\mathrm m}$ 
or $\mathbb{G}_{\mathrm a}$, 
or in other words, $\Aut(X)$ contains 
a non-torsion algebraic element.
\end{conj}
Due to Theorem \ref{thm:PZ} below, 
Conjecture \ref{conj-6} 
holds for affine surfaces. It holds also provided 
$\Aut(X)$ contains a  nontrivial connected set of 
commuting elements. Indeed,  we have 
the following results. 
\begin{thm}[{\rm\cite[Theorem B]{CRX19}}]
\label{thmB-comm}
Let $Y$ be an irreducible algebraic subvariety 
of $\Aut(X)$ such that $Y$ contains the identity
and any $g_1,g_2\in Y$ commute. 
Then the subgroup of $\Aut(X)$ generated 
by $Y$ is an abelian connected algebraic 
subgroup of $\Aut(X)$.
\end{thm}
\begin{cor}  Let $G$ be a connected 
abelian subgroup 
of $\Aut(X)$.
\begin{itemize} 
\item If $G$ is nontrivial, then it contains 
a subgroup isomorphic to either 
$\mathbb{G}_{\mathrm m}$ or 
$\mathbb{G}_{\mathrm a}$.
\item If $G$ is a closed ind-subgroup, 
then $G$ is nested. 
\end{itemize}
\end{cor}
As a consequence, an eventual counterexample 
to Conjecture \ref{conj-6} 
has the property that every commutative subgroup 
of $\Aut(X)$ is discrete, 
or, in other terms, every algebraic element  of 
$\Aut(X)$ is a torsion element. 
See also \cite[Theorem A]{RvS21} for a description 
of maximal abelian subgroups of $\Aut(X)$ which 
consist of unipotent elements. 
\begin{que}\label{ques-6}
Assume that the group $\Aut(X)$ 
does not contain algebraic 
elements different from $\id_X$.
 Is it then discrete?
\end{que}
Notice that by Theorem \ref{thmB-comm}, 
 if  $\Aut(X)$ has no 
 algebraic element different from $\id_X$, then
every abelian subgroup $H$
of $\Aut(X)$ has no torsion 
and the neutral component $H^0$ is trivial. 
Cf. also \cite[Problem 20]{Arz18a}.
\section{Nilpotency and solvability}
\begin{thm}[{\rm\cite[Theorem B]{KrZa24}}]
Let $G \subset \Aut(X)$ be an algebraically 
generated subgroup, see subsection \ref{ss:ags}. 
If $G$ is solvable, then the following hold.
\begin{itemize}
\item[(a)] $G=U\rtimes T$ where $T$ 
is an algebraic torus and $U=R_{\mathrm u}(G)$ 
is a nested unipotent group. 
\item[(b)] $G$ is unipotent if and only if the 
generating subgroups are. 
\item[(c)] If, in addition, $G$ is generated 
by a finite family of connected algebraic groups, then 
$G$ is a connected algebraic group. 
\end{itemize}
\end{thm}
\begin{thm}[{\rm\cite[Theorem C]{KrZa24}}]
A nested unipotent subgroup $U \subset \Aut(X)$ is 
solvable of derived length 
$\leq \max\{\dim Ux\mid x \in X\}\leq \dim X$.
\end{thm}
\begin{rem}\label{rem2} 
A unipotent algebraic group is nilpotent.
By contrast, a nested unipotent group $G=\cup_i G_i$
 is not necessarily nilpotent, since
the nilpotency class of the unipotent algebraic subgroups
$G_i$ might not be bounded.

For example, consider  the 
\emph{de Jonqui\`eres subgroup} 
 $\Jonq_n\subset\Aut(\AA^n)$ of 
triangular automorphisms of the form
\[
\phi=(a_1x_1+f_1, a_2x_2 +f_2(x_1),
\ldots, 
a_n x_n + f_n(x_1,\ldots,x_{n-1}))
\]
where $a_i \in \kk^*$ and 
$f_i\in\kk[x_1,\ldots,x_{i-1}]$. 
The unipotent radical 
$\JJJ_n=R_{\mathrm u}(\Jonq_n)$ 
is a closed nested unipotent subgroup 
of $\Aut(\AA^n)$. 
For $n\ge 2$ its Lie algebra $\Lie(\JJJ_n)$ 
is not nilpotent. 
However, the Lie algebra of a nilpotent ind-group 
is nilpotent, see \cite[Lemma 5.1.4(3)]{KrZa24}. 

Indeed, let for instance $n=2$ and
\[L_d = \langle \partial/\partial x_1, 
x_1^d\partial/\partial 
x_2\rangle_{Lie}\subset \Lie(\JJJ_2).\]
It is easily seen that the $d$th member $(L_d)_d$ of 
the lower central series of $L_d$ does not vanish, 
and so $(\Lie(\JJJ_2))_d\neq 0$ for any $d\geq 1$, 
 cf. \cite[Remark 5.3.3]{KrZa24}. 
\end{rem}
\begin{que} 
Let $\G$ be a connected ind-group. 
Is it true that $\G$ is solvable (resp., nilpotent) 
if $\Lie\G$ is? 
\end{que}
\begin{defn} Let us say that a subgroup 
$G\subset\Aut(X)$ is 
\emph{a-generated} if $G=\langle Y\rangle$ 
where $Y\subset\Aut(X)$ 
is an irreducible algebraic subset containing $\id_X$. 
Replacing $Y$ by $Y\cdot Y^{-1}$ we may assume 
that $Y$ is symmetric, that is, $Y=Y^{-1}$.
\end{defn}
\begin{conj}\label{conj:a-generated} If $G$ 
is a-generated 
and nilpotent (resp., solvable) 
then $G$ is algebraic. 
\end{conj}
In the case of a nilpotent (resp., a solvable) 
subgroup $G$ one can try to proceed by induction 
on the nilpotency class of $G$ 
(resp., on the derived length 
of $G$), where the case $n(G)=1$
(i.e. $G$ abelian) is settled by 
Theorem \ref{thmB-comm}.
In the general case, the induction works
provided the following conjecture is true. 
\begin{conj} If $G$ is a-generated and 
nilpotent (resp., solvable) 
then any member of its lower central series 
(resp., derived series) is a-generated as well. 
\end{conj}
Let $z(G)$ stands for the center of $G$.
If $G/z(G)$ is abelian, then Conjecture 
\ref{conj:a-generated} 
holds if the following is true:
\begin{conj} Consider an admissible filtration 
$\Aut^0(X)=\bigcup A_n$ 
by irreducible affine subvarieties $A_n$. Let $G$ be 
a (nilpotent, resp., solvable) subgroup of $\Aut^0(X)$ 
such that $G=\langle G\cap A_n\rangle$ 
for some $n$. 
Then $z(G)=\langle z(G)\cap A_m\rangle$ 
for some $m$ 
(resp., $(G,G)=\langle (G,G)\cap A_m\rangle$ 
for some $m$).
\end{conj}
 The following question concerns
a modified version of stable triangulation, 
see Section \ref{ss:triangulation}. 
\begin{que}[\rm{\cite[Sect.~5.2, Question 6]{KrZa24}}]
\label{triangulation}  
Let $G \subset \Aut(X)$ 
be a solvable or nilpotent connected subgroup. 
Does there exists a closed embedding $X \into \AA^n$ 
such that $G$ extends to a subgroup of the de Jonqui\`eres 
group $\Jonq_n$? Is this true if $G$ is 
 algebraic or nested?
\end{que}
 \begin{rem} \label{rem-6.10}
The answer to the latter question is affirmative 
in a particular  setup where $X=\AA^k$ 
and $G\subset\Aut(\AA^k)$ 
is a unipotent algebraic group, see 
\cite[Remark 1.4.2(3)]{KrZa24}. For instance, 
the Nagata-Bass $\Ga$-subgroup 
of $\Aut(\AA^3)$ can be triangulated in $\AA^4$ 
in the above sense, being non-stably 
triangulable in the usual sense, 
cf. Example \ref{ex:NB}. 
\end{rem}
Denote by $\LND(X)\subset\VF(X)$ the set of all 
locally nilpotent derivations of $\OOO(X)$. 
For a vector space $V$ over $\kk$ we let $\LND(V)$ 
be the set of locally nilpotent endomorphisms of $V$. 
\begin{defn}[{\cite{Dai19}}] 
A subset $Z$ of endomorphisms of a vector space 
$V$ is called 
\emph{uniformly locally nilpotent},  or  ULN
for short, if for any $v\in V$ 
there is $n\in\NN$ such that $a_1\cdots a_n(v)=0$ 
whatever are $a_1,…,a_n \in Z$.  
\end{defn}
Clearly, if $Z$ is ULN then  $Z\subset \LND(V)$ and 
${\rm span}_{\kk}\,(Z)\subset\End_\kk(V)$ is ULN too. 

Let $A\subset\Aut(X)$ be  an irreducible algebraic subvariety 
which contains $\id_X$,
and let $\langle A\rangle$ be the subgroup of 
$\Aut(X)$ generated by $A$. 
The notation
$L_A\subset\VF(X)$ 
 has the same meaning as in subsection \ref{ss:ags}. 
Given an ind-subgroup $\G$, 
 the \emph{adjoint representation} 
 $\Ad(\G)\colon\G\to \GL(\Lie\G)$
 is the tangent representation 
 at the fixed point $\id_X\in\G$
of the action of $\G$ on itself by inner automorphisms
 ${\rm Int}(g)\colon h\mapsto g\cdot h\cdot g^{-1}$, 
 see \cite[Section 7.5]{FuKr18}.
\begin{que} 
Do the following equivalences hold?
\begin{enumerate}
\item
$L_A\subset \LND(X)
\Leftrightarrow\Ad(A)\subset \End_{Lie}(\VF(X))$ 
is {\rm ULN} 
$\Leftrightarrow\langle A\rangle$ 
is a unipotent algebraic group;
\item  $\langle A\rangle$ is metabelian 
$\Leftrightarrow\langle\Ad(A)\rangle_{Lie}$ 
is commutative; 
\item $\langle A\rangle$ is solvable 
$\Leftrightarrow\langle\Ad(A)\rangle_{Lie}$ is solvable;
\item $\langle A\rangle$ is nilpotent 
$\Leftrightarrow\langle\Ad(A)\rangle_{Lie}$ is nilpotent.
\end{enumerate}
\end{que}
\section{$\Aut$-quasihomogeneous and flexible
affine varieties}
\subsection{$\Aut$-quasihomogeneous affine varieties}\label{ss:quasi-homogeneous}
\begin{defn} We say that a 
(not necessarily affine) variety $X$ is 
\emph{$\Aut$-homogeneous} 
(resp. \emph{$\Aut$-quasihomogeneous}) if $\Aut(X)$ 
acts transitively on $X$ (resp. acts on $X$ 
with an open orbit). 
\end{defn}
The following characterizations of 
$\Aut$-quasihomogeneity for smooth  affine surfaces 
are due to Gizatullin \cite{Giz71b}.
This was extended to normal affine surfaces by 
Dubouloz \cite{Dub04}. 
\begin{thm}
Let $X$ be a normal affine surface $X$ 
non-isomorphic to 
$\AA^1 \times (\AA^1\setminus\{0\})$. 
Then the following are equivalent:
\begin{itemize}
\item $X$ is $\Aut$-quasihomogeneous; 
\item there are two effective $\Ga$-actions on $X$ 
with distinct general orbits;
\item $X\cong \bar X\setminus D$ where $\bar X$ 
is a normal projective surface and 
$D\subset \bar X$ 
is a reduced divisor 
with simple normal crossings 
contained in the smooth locus 
$ {\rm reg}(\bar X)$ of $\bar X$, 
which consists of
a chain of smooth rational curves. 
\end{itemize}
\end{thm}
Notice that many, but not all, smooth Gizatullin 
surfaces are $\Aut$-homogeneous. 
See
Gizatullin \cite{Giz71a} and Popov \cite{Pop73}
for a list of normal affine surfaces that 
are (quasi)homogeneous spaces of algebraic groups;
see also \cite[Theorem 4.8]{KPZ17}. 

In higher dimensions there is a classification of 
$\Aut$-homogeneous affine  toric varieties 
and complete toric varieties, 
see \cite{Arz18b};
cf. also the recent preprint \cite{ArZa23}.
\subsection{Flexible affine varieties}
Flexible affine varieties were introduced in 
\cite{AFK$^+$13a},
see also \cite{KZ99} and \cite{AKZ12}. 
For overviews on the subject
and examples we address \cite{Arz23}, 
\cite{AFK$^+$13b}, 
 and \cite[Section 7]{CPPZ21}.
\begin{defn} Let $X$ be a  quasiaffine variety. 
We say that a smooth point $x\in X$ is \emph{flexible} 
if the tangent space $T_xX$ 
is spanned by the tangent vectors to the orbits $Ux$ 
of $\Ga$-subgroups $U$ of $\Aut(X)$. 
The variety $X$ is called \emph{flexible} 
if every smooth point $x\in X$ is flexible. 
Clearly, $X$ is flexible if
at least one smooth point $x\in X$  
is and the group $\Aut(X)$ acts 
transitively on the smooth locus $ {\rm reg}(X)$.
\end{defn}
We let $\SAut(X)\subset\Aut^0(X)$ be the 
 subgroup generated by all 
 $\Ga$-subgroups of $\Aut(X)$, 
in other words, by all unipotent elements of $\Aut(X)$ 
\footnote{The subgroup 
$\SAut(X)$ is denoted by ${\rm SAut}(X)$ 
in \cite{AFK$^+$13a}.
We prefer the notation $\SAut(X)$ 
in order to keep ${\rm SAut}(\AA^n)$ 
for the subgroup of $\Aut(\AA^n)$ of 
volume-preserving automorphisms.}. 
We also let $\Aut_{\rm alg}(X)$ be the subgroup 
generated by all algebraic elements of $\Aut(X)$. 
Notice that $\Aut_{\rm alg}(X)$ is not 
necessarily connected,
and every algebraic element
$a\in \Aut_{\rm alg}(X)\setminus 
\Aut^0(X)$ is a torsion element. 

The following theorem (conjectured 
in \cite[Sect. 4.2]{AKZ12}) 
was proven in \cite{AFK$^+$13a} for affine varieties and 
 extended in \cite[Theorem 2]{APS14} and 
 \cite[Theorem 2.12]{FKZ16} to quasiaffine varieties. 
\begin{thm}\label{thm:flex}
For a quasiaffine variety $X$ of dimension at least 2, 
the following are equivalent:
\begin{itemize}
\item $X$ is flexible;
\item the group $\SAut(X)$ acts transitively 
on the smooth locus ${\rm reg}(X)$; 
\item the group $\SAut(X)$ acts infinitely 
transitively on ${\rm reg}(X)$. 
\end{itemize}
\end{thm}
The flexibility survives when removing a subvariety 
of codimension $\ge 2$. 
\begin{thm}[{\rm\cite[Theorem 0.1]{FKZ16}}] 
Let $X$ be a flexible smooth quasiaffine variety 
of dimension $\ge 2$, and let $Y\subset X$ 
be a closed subset of codimension $\ge 2$. 
Then $X\setminus Y$ is flexible.  
\end{thm}
As an application of Theorem \ref{thm:flex}, 
let us mention the following two interpolation results. 
The  first follows immediately from 
\cite[Theorem 4.14 and Remark 4.15]{AFK$^+$13a}.
\begin{thm}\label{thm:interp}
Let $X$ be a flexible affine variety of 
dimension $n\ge 2$. 
Then for every $m \ge 0$ and every finite subset 
$Z  \subset {\rm reg}(X)$ there exists 
an automorphism 
$g\in \SAut(X)$ with prescribed $m$-jets 
at the points 
$p \in Z$ provided each jet fixes its center $p$.
\end{thm}
\begin{cor}
Under the assumptions of Theorem \ref{thm:interp} 
there exists an $\AA^1$-orbit of a 
$\Ga$-action on $X$ 
which passes through every point $p\in Z$ 
and interpolates a given smooth curve jet at $p$. 
\end{cor}
See also \cite[Theorem 1.15]{AFK$^+$13a} 
for an analogue
of Kleiman's Transversality theorem for 
flexible varieties. 

The first nontrivial examples of flexible varieties 
appeared in \cite{Giz71a} and \cite{Giz71b}; 
these are $\Aut$-homogeneous Gizatullin affine surfaces, 
see subsection \ref{ss:quasi-homogeneous}.
See \cite{GD75}, \cite{GD77}, \cite{KPZ17} and references therein
for systematic studies on the automorphism groups 
of Gizatullin surfaces. 

Let us mention two interesting examples 
of flexible affine varieties of the form
$X_{n,d}=\PP^n\setminus H_d$, where $H_d$  is a hypersurface
of degree $d$ in $\PP^n$ with
$(n,d)\in\{(26,3),\,(55,4)\}$, see
\cite[Proposition 3.7 and Theorem 4.3]{Giz18}. 
These are homogeneous spaces of simple linear groups 
of type $E_6$ and $E_7$, respectively. The cubic form that defines 
$H_3\subset\PP^{26}$
was first written by Camille Jordan in 1870. It is related to
the configuration of $27$ lines on a smooth 
cubic surface in $\PP^3$. 
The quartic form  that defines $H_4\subset\PP^{55}$ is related to
the $56$ exceptional curves on a Del Pezzo surface of degree 2, 
and to the 28 bitangent lines of a plane quartic curve.
\subsection{Flexibility of toric varieties}
Let $\TT=(\Gm)^n$ stand for an 
algebraic $n$-torus.
Recall that a $\TT$-variety $X$ of dimension $n$ is 
a \emph{toric variety} if $\TT$ 
acts effectively on $X$  with an open orbit.  
One says that a toric variety $X$ 
has a torus factor 
if $X\cong Y\times (\Gm)^k$ 
for some toric variety $Y$ and some $k>0$,
or equivalently, if there is a nonconstant 
invertible regular function on $X$.

Any toric affine variety $X$  of dimension $n\ge 2$  
with no torus factor
can be defined by a strictly convex polyhedral cone 
in the integer lattice $\ZZ^n$. For instance, 
the affine plane $\AA^2$ 
with the standard diagonal torus action is defined 
by the first quadrant $(\NN\cup\{0\})^2$
of $\ZZ^2$.
\begin{defn}\label{Demazure} Let $X$ be a toric variety. 
A \emph{Demazure root subgroup} of $\Aut(X)$ 
is a $\Ga$-subgroup $U\subset\Aut(X)$
which is normalized by the acting torus $\TT$. 
\end{defn}
Such subgroups, or rather the associated 
locally nilpotent derivations of $\OOO(X)$,
correspond to certain lattice points of $\ZZ^n$ called 
\emph{Demazure roots}. In particular, 
if $X=\AA^2$ then the Demazure roots are 
the lattice points $(-1,k)$ and $(l,-1)$ 
with $k,l\ge 0$. They correspond to 
the locally nilpotent 
derivations $y^k\p/\p x$ and $x^l\p/ \p y$, 
respectively. 
It is known that any normal toric affine variety 
of dimension  $\ge 2$ 
with no torus factor is flexible, 
see \cite{AKZ12}. The following theorem says more. 
\begin{thm}[{\rm \cite[Theorem 1.3]{AKZ19}}]
\label{thm:main-2}
Let $X$ be a toric affine variety of dimension 
at least 2. 
If $X$ has no torus factor and is smooth in 
codimension $2$, 
then one can find 
a finite collection of Demazure root subgroups 
such that the group generated 
by these acts infinitely transitively 
on the smooth locus 
${\rm reg}(X)$. 
\end{thm}
\section{Rigid affine varieties} 
\begin{defn} 
An affine variety $X$ is called \emph{rigid} 
if it admits no effective $\Ga$-action.
\end{defn}
\begin{exa} Conjecturally, a Pham--Brieskorn hypersurface 
$X\subset \AA^{n+1}$, $n\ge 2$,
defined by
\[x_0^{a_0} +x_1^{a_1} +\cdots +x_n^{a_n} =0
\quad\text{where}\quad
2\le a_0\le a_1\le\cdots\le a_n\] 
is rigid if and only if $a_1\ge 3$, see 
\cite[Conjecture~1.22]{CPPZ21}.
This is indeed the case for $n=2$ 
\cite[Lemma~4]{KZ00} and for $n=3$
\cite[Main theorem]{ChDu23}. 
For example, the Fermat threefold hypersurface
\[x_0^3+x_1^3+x_2^3+x_3^3=0\] in $\AA^4$ is rigid;
see \cite[Theorem~1.13]{CPPZ21} for a far more general result. 
See also \cite{EGS23} and references therein for 
a description of rigid trinomial affine varieties.  
\end{exa}
For the following conjecture see e.g. 
\cite[Conjectures 1.0.1 and 1.0.2]{PZ22};  cf. also 
Conjecture \ref{conj-6} and Question \ref{ques-6}.
\begin{conj} \label{conjPZ}
Let $X$ be a rigid affine variety. Then $\Aut^0(X)$ 
is an algebraic torus (of rank $\le\dim(X)$).
\end{conj}
Conjecture \ref{conjPZ} is confirmed
for rigid toric varieties and varieties with 
a torus action of complexity one, see \cite{BG23}.
Due to the following theorem, this conjecture also holds
in dimension $2$. 
\begin{thm}[{\rm \cite[Theorem 1.3(1)]{PZ22}}]  
\label{thm:PZ}
Assume that the base field $\kk$ is uncountable. 
 Let $X$ be a normal affine surface over $\kk$. Then
$\Aut^0(X)$ is an algebraic group if and only 
if $X$ is rigid, 
if and only if $\Aut^0(X)$ is an algebraic torus 
(of rank $\le 2$).
\end{thm}
We have the following geometric criterion of rigidity 
of a normal affine surface.
\begin{thm}[{\rm \cite[Theorem 1.3(2)]{PZ22}}] 
Assume that the base field $\kk$ is uncountable. 
Then a normal affine surface $X$ over $\kk$ is rigid 
if and only if 
$X$ admits a relatively minimal completion $\bar X$ 
by a reduced effective normal crossing divisor 
$D$ supported on 
${\rm reg}(\bar X)$ such that 
any extremal chain of the dual graph 
$\Gamma(D)$ which consists 
of rational components $C_1,\ldots,C_k$ 
satisfies $C_i^2\le -2$ for $i=1,\ldots,k$. 
\end{thm}
The dual graph $\Gamma(D)$ is defined as follows: 
its vertices  correspond 
to components of $D$
and the edges correspond to the pairs of 
components that intersect. A chain $L$ in 
$\Gamma(D)$ is called \emph{extremal} if 
all vertices of $L$ are of valence $\le 2$ 
in $\Gamma(D)$
and at least one vertex is of valence 1. 
\begin{exa}[{\rm \cite[Example 4.3.6]{PZ22}}] 
Consider the affine surface 
$X=\mathbb{P}^2\setminus {\rm supp}(D)$ 
where $D$ is 
a reduced effective divisor on $\mathbb{P}^2$ 
with only nodes as singularities. 
Then $X$ is not rigid if and only if $\deg(D)\le 2$, 
that is, $D$ is either 
a projective line of a (possibly, singular) conic.
\end{exa}
By definition, a nested ind-group consists of 
algebraic elements. 
The following theorem provides a partial converse.
\begin{thm}[{\rm\cite[Theorem 1.1]{PR24}, 
cf. also \cite{PR23}}]
\label{thm:PR22}
For an affine variety $X$ the following 
conditions are equivalent:
\begin{itemize}
\item $\Aut^0(X)$ consists of algebraic elements;
\item $\Aut^0(X)$ is a closed nested ind-subgroup 
of $\Aut(X)$;
\item $\Aut^0(X) = \SAut(X)\rtimes T$ 
where $T$ is a maximal torus of $\Aut(X)$ and 
$\SAut(X)=R_u(\Aut(X))$ is abelian
and consists of all unipotent elements of $\Aut(X)$.
\end{itemize}
\end{thm}
\begin{cor}[{\rm \cite[Corollary 4.3]{PR24}}]
\label{cor:PR24}
Let $X$ be a rigid affine variety. If $\Aut^0(X)$ 
consists of 
algebraic elements, then $\Aut^0(X)$ is 
an algebraic torus 
of dimension $\le \dim(X)$. 
\end{cor}
The proof exploits the following result of 
independent interest.
\begin{thm}[{\rm \cite[Theorem 3.1]{PR24}}]
\label{thm2:PR22}
Let $\G$ be a connected ind-group and 
$\H \subset \G$ be 
a closed connected nested ind-subgroup. 
Assume that any $g \in \G$ is a torsion 
element modulo $\H$, 
that is, $g^d \in \H$ 
for some $d=d(g)>0$. Then $\H=\G$.
\end{thm}
The following theorem gives 
a necessary condition for rigidity. 
For the proof see \cite[Theorem 2.1]{AG17}; 
cf. also \cite[Theorem 3.3]{FZ05b} 
in the case of affine surfaces. 
\begin{thm}
Suppose that $\Aut(X)$ contains
 two non-commuting $\Gm$-subgroups. 
Then $\Aut(X)$ contains a $\Ga$-subgroup. 
\end{thm}
This theorem implies (a) below;
see  \cite[Proposition 5.1 
and Remark 5.2]{PR24} for the proof of (b).
\begin{thm}\label{thm:AG}
Let $X$ be a rigid affine variety. Then  the following hold.
\begin{itemize} 
\item[(a)] There exists an algebraic torus 
$\TT\cong(\Gm)^k$
 in $\Aut(X)$
which contains any other algebraic subtorus
of $\Aut(X)$. 
In particular, $\TT$ is a normal subgroup 
of $\Aut(X)$.
\item[(b)] $\TT$ 
is a central subgroup of $\Aut^0(X)$ and 
the only maximal connected abelian 
ind-subgroup of $\Aut^0(X)$.
\end{itemize}
\end{thm}
The next theorem 
confirms Conjecture 
\ref{conjPZ} in a particular setting of rigid varieties 
with a torus action. 
For statement (a) see
\cite[Example 2.3 and Remark 2.4]{AG17} 
and \cite[Sect. 2]{AKZ12};  (b) is proven in
 \cite[Theorem 6.4]{BG23}.
Recall that the \emph{complexity} 
of an algebraic group action 
is the codimension of general orbits. 
\begin{thm}\label{thm:BG}
\begin{itemize} 
\item[(a)] Let $X$ be a normal toric variety. 
If $X$ is rigid, 
then $X\cong\TT$ is a torus and 
$\Aut^0(X)\cong\TT$ 
is the torus acting naturally on itself.
\item[(b)]
Let $X$ be a  normal rational
affine variety with a torus action 
of complexity one
and a finitely generated divisor class group.
Suppose that $X$ is rigid 
and every invertible regular function 
on $X$ is constant. 
Then $\Aut^0(X)$ 
is an algebraic torus.
\end{itemize}
\end{thm}
Concerning the conjugacy of tori in $\Aut(X)$ 
we have the following general results, see e.g. 
\cite[Proposition 10.5.4]{FuKr18} for (a)
and \cite{Dem82}, \cite{Gub98} and \cite{BH03}
for (b).
\begin{thm} $\,$
\begin{itemize}
\item[(a)] Let $T \subset \Aut(X)$ be an 
algebraic torus. 
Then $\dim(T) \le \dim(X)$. 
\item[(b)]
Assume that $X$ is a normal toric variety.
Then all algebraic tori $T\subset \Aut(X)$ 
of dimension $\dim(X)$ are conjugate. 
Furthermore,
any algebraic torus of dimension 
$\dim(X)-1$ in $\Aut(X)$ is contained in  
an algebraic torus of dimension $\dim(X)$.
\end{itemize}
\end{thm}
However,
for a non-toric affine variety $X$ the group $\Aut(X)$
can contain several pairwise non-conjugate 
maximal algebraic tori. 
\begin{exa}\label{ex:Dani} $\,$
1.  (See \cite[Theorem 2]{Dani89}.)
 Let $D$ be the Danielewski surface 
$\{xy-z^2+z=0\}$ in $\AA^3_{\CC}$
and let $X=D\times\AA^1$. 
Then $\Aut(X)$ contains an infinite sequence of pairwise 
non-conjugate algebraic $2$-tori 
and no $3$-torus.
Furthermore, let $X$ be the underlying affine variety of 
$\SL(2,\CC)\times\AA^1_{\CC}$. 
Then the number of conjugate classes of maximal
tori in $\Aut(X)$ is infinite \cite[P.~7]{Dani89}.

2. There exist examples of $\Aut$-homogeneous 
smooth non-toric affine surfaces $X$
such that $\Aut(X)$ contains an algebraic 
family of pairwise 
non-conjugate algebraic $1$-tori depending 
on 1 or 2 parameters, 
see \cite[Theorem 1.0.1]{FKZ11} and also 
\cite[Theorem 5.5]{KPZ17}. 
Cf.\ also Example \ref{ex:DG} below of 
Danilov-Gizatullin surfaces,
where the number of conjugate classes 
of $1$-tori  is finite.
\end{exa}
\begin{rem}[{\rm 
see \cite{KML97} and \cite{DMJP10}}] 
Recall 
that the Koras-Russel threefold
$X$
given in $\AA^4_{\CC}$
by equation $x+x^2y+z^2+t^3=0$ is
diffeomorphic to $\AA^3_{\CC}$, 
but non-isomorphic to $\AA^3_{\CC}$. 
It is neither flexible, nor rigid. The
 group $\Aut(X)$ is infinite dimensional, 
 acts on $X$ with an open orbit 
 and a unique fixed point (the origin). It
 coincides with the stabilizer of $X$ in $\Aut(\AA^4_{\CC})$. 
This group is naturally isomorphic 
to the subgroup of $\Aut(\AA^3_{\CC})$
 (where $\AA^3_{\CC}=\Spec\CC[x,z,t]$)
 of all automorphisms which leave invariant 
 the ideals $(x)$ and 
$(x^2, x+z^2+t^3)$.
\end{rem}
\section{Chapters on $\Aut(\AA^n)$}
\subsection{Groups acting infinitely 
transitively on $\AA^n$}
Recall that the \emph{root subgroups} 
of $\Aut(\AA^n)$
are the $\Ga$-subgroups normalized 
by the standard $n$-torus, see Definition \ref{Demazure}. 
\begin{thm}[{\rm \cite{And19}, \cite{AKZ19}}]
\label{thm:main-1} 
For any $n\ge 2$ one can find three $\Ga$-subgroups 
(resp., $n+2$ root subgroups) of $\Aut(\AA^n)$ 
which generate a subgroup  acting  
infinitely transitively on $\AA^n$. 
\end{thm}
\subsection{Tame subgroups}\label{ss:tame}
Notice that the de Jonqui\`eres subgroup 
$\Jonq_n\subset\Aut(\AA^n)$ (see Remark 
\ref{rem2})
 is generated by the standard torus 
 $\TT$ of $\GL(n,\kk)$
and the unipotent radical 
$\JJJ_n=R_{\mathrm{u}}(\Jonq_n)$.
In turn,  $\JJJ_n$ is generated by the 
 triangular root subgroups
\[\exp(tx_{1}^{m_{1}}\cdots x_{i-1}^{m_{i-1}}\p_i),
\quad i=1,\ldots,n\quad\text{where}\quad
\p_i=\p /\p x_i.\] 
The Lie algebra 
\[\Lie(\JJJ_n)=\kk\p_1\oplus \kk[x_1]\p_2\oplus
\cdots\oplus\kk[x_1,\ldots,x_{n-1}]\p_n\]
consists of locally nilpotent derivations. 
This is a maximal Lie subalgebra  in $\VF(\AA^n)$ 
consisting of  locally nilpotent derivations, 
see \cite[Theorems 1 and 7]{Sku21};
cf. also \cite[Theorems 1.4 and 1.5]{BPZ24} 
and the references therein for closely related results.
For $\kk=\CC$, $\Jonq_n$ 
is a Borel subgroup of $\Aut(\AA_{\CC}^n)$, that is, 
a maximal (connected) solvable subgroup, 
see \cite[Corollary 1.2]{FuPo18}.
\begin{defn}
An automorphism $\phi\in\Aut(\AA^n)$ 
is said to be \emph{tame} if it
is a composition of affine and triangular 
automorphisms; 
otherwise it is called \emph{wild}. 
The tame automorphisms form  the 
\emph{tame subgroup} 
$\Tame_n=\langle \Aff_n, \Jonq_n\rangle$ 
of $\Aut(\AA^n)$.
\end{defn}
For the following classical theorem 
see \cite{Jun42}, \cite{vdK53}; cf.  also 
\cite[Theorem 2]{Kam75} and 
\cite[Theorem 2.3.8]{Miy24}.
\begin{thm}[Jung-van der Kulk]\label{thm:JvdK} 
Let $k$ be an arbitrary field.
Then $\Aut(\AA^2_{k})=\Tame_2$.
Moreover, $\Aut(\AA^2_k)$ is a free product  
$\Aff_2*_C\Jonq_2$ amalgamated 
over $C=\Aff_2\cap \Jonq_2$. 
\end{thm} 
Using this amalgam, Danilov \cite{Dan74} 
showed that the 
subgroup 
${\rm SAut}(\AA^2_\kk)=\{f\in\Aut(\AA^2_\kk)|\det(f)=1\}$
of unimodular automorphisms is not a simple group. 
 
The things are different for $\Aut(\AA^3_\kk)$.  
Nagata considered in \cite[Section 2.1]{Nag72} 
the following automorphism of 
$\AA^3_\kk=\Spec(\kk[x,y,z])$:
 \[\Phi\colon (x,y,z)\mapsto 
(x+\sigma z, y+2\sigma  x+\sigma^2 z, z)
\quad\text{where}\quad\sigma=x^2-yz.\]
We have 
\[\Phi=\exp(\sigma\p)\quad\text{where}
\quad \p=z\p/\p x+2x\p/\p y.\]
Indeed, $\p(\sigma)=0$ and so, 
$\Phi^*(\sigma)=\sigma$.
Since also $\Phi^*(z)=z$ we have 
$\Phi\in\Aut_{\kk[z]} \kk[z][x,y]$. 
 Nagata showed that $\Phi$ is not generated 
by the affine and triangular automorphisms that fix $z$,
and suggested that $\Phi$ is wild, i.e., 
$\Phi\notin \Tame_3$, 
see \cite[Conjecture 3.1]{Nag72}. 
The following theorem due to Shestakov 
and Umirbaev confirms Nagata's Conjecture.
\begin{thm}[{\rm \cite{SU04a, SU04b}; 
see also \cite[Chapter 1]{vdEKC21}}]
The Nagata automorphism  $\Phi$
is wild. So, $\Tame_3$ is a proper subgroup 
of $\Aut(\AA^3)$. 
\end{thm}
Notice that $\Tame_3$ is not a closed 
ind-subgroup of $\Aut(\AA^3)$, 
see \cite{EP15}. It is unknown whether 
$\Tame_3$ is dense in $\Aut(\AA^3)$. 
\begin{exa}[{\rm \cite[Sect.~17.3]{FuKr18}}]
\label{exa:FK18}
Letting $\AA^3=\Spec\kk[x,y,z]$ consider  
the closed connected ind-subgroup 
$\G\subset\Aut(\AA^3)$ consisting 
of the automorphisms 
which send $z$ to $z$.
Then $\G^t=\G\cap\Tame_3$ is 
a proper closed ind-subgroup 
of $\G$. Nevertheless, their Lie algebras 
coincide:
$\Lie(\G^t)=\Lie(\G)$, contrary to the claim 
of \cite[Theorem 1]{Sha81}. 
\end{exa}
\begin{que} 
Do we have $L_{\G}=\Lie \G$ and 
$L_{\G^t}=\Lie\G^t$?
Is it true that  $L_{\G}=L_{\G^t}$?
\end{que}
Let us mention also the following result.
\begin{thm}[{\rm\cite[Main Theorem]{KrSt13}}]
\label{thm:inner}
Every automorphism of the group $\Aut(\AA^n_\CC)$
is inner up to field automorphisms when restricted to 
the tame subgroup $\Tame_n$.
\end{thm}
\begin{rems}\label{rem:Anick} 
1. For $n=2$ the assertion 
was established in \cite{Des06}.
Since $\Tame_2=\Aut(\AA^2_\CC)$,
see Theorem \ref{thm:JvdK}, it follows that
every automorphism of 
$\Aut(\AA^2_\CC)$ is inner 
up to a field automorphism. 
The latter does not hold, in general, 
for toric affine surfaces. For example, 
this fails for the quotient of $\AA^2$ 
by the cyclic group generated by
$g\colon (x,y) \mapsto (\xi^e x,\xi y)$, 
where $\xi$ is a $d$th primitive root of unity, 
$1 < e < d$, $(e, d) = 1$ and 
$e^2 \not\equiv 1 \mod d$, see 
\cite[Remark 5.13]{LRU23}.

2. Any automorphism of $\Aut(\AA^n)$ 
as an ind-group is inner, 
see \cite{BKYE16}. 
The same holds for the automorphisms 
of the subgroup ${\rm SAut}(\AA^n)$ 
of volume-preserving automorphisms of $\AA^n$,
see \cite{Kra17}.  See also \cite{Sta13}, 
\cite[Theorem 1.3]{UZ21} 
and the literature therein for some generalizations 
of Theorem \ref{thm:inner}.

3. A tame subgroup was defined for 
some other automorphism groups,
for instance,
\begin{itemize} 
\item for $\Aut(Q_3)$ where $Q_3$ 
is a smooth 
affine quadric threefold
realized as the underlying variety of 
$\SL(2,\CC)$, 
see \cite{LV13}; 
\item for $\Aut(A_n)$ where $A_n=
\kk\langle x_1,\ldots,x_n\rangle$ 
is the free associative algebra in $n\ge 2$ 
variables,
see \cite{BKY12b};
 \item for $\Aut(X)$, 
where $X$ is a toric affine variety, 
see \cite[Definition 5.3]{AG10}. 
\end{itemize}
A theorem due to Makar-Limanov and 
Czerniakiewicz says
that $\Aut(A_2)=\Tame(\Aut(A_2))$, 
see \cite{ML70}, \cite{Cze71-72} and 
\cite[Theorem~4.1]{Coh85}.
At the same time,
$\Tame(\Aut(Q_3))$ and $\Tame(\Aut(A_3))$ 
are proper subgroups
of $\Aut(Q_3)$ and $\Aut(A_3)$, respectively,
see \cite[Theorem 6.1]{AG10} 
and \cite[Section 5.1]{LV13} for the former and 
\cite[Corollary 2.1]{Umi07} for the latter.
The role of the Nagata automorphism for 
$\Aut(A_3)$
is played by
the Anick automorphism
\[(x,y,z)\mapsto (x+z\sigma,y+x\sigma,z)
\quad\text{where}\quad \sigma=xz-zy,\]
which is wild. See also \cite{SZ24} for 
a wild (but stably tame) automorphism of Anick type
of a $3$-generated free Poisson algebra, 
which induces a tame automorphism 
of the 3-generated polynomial algebra.  

Likewise, consider the quadratic cone 
$S_3\subset\AA^4$ 
given by equation 
$x_1x_4 - x_2x_3=0$. 
There exists an automorphism $\tau$ 
of $\AA^4$ 
which restricts to $S_3$
and lifts to the Cox ring $\kk[x_1, x_2, x_3, x_4]$ 
of $S_3$
yielding the Anick automorphism of $\AA^4$:
\[(x_1, x_2, x_3, x_4) \mapsto (x_1, x_2 + x_1\delta, 
x_3, x_4 + x_3\delta) \quad\text{where}
\quad \delta=x_1x_4 - x_2x_3.\]
The restriction $\tau|_{S_3}$ is a wild automorphism 
of $S_3$, see \cite[Theorem 6.1]{AG10}.
\end{rems}
The following theorem relates 
the automorphism group of $\AA^n$ 
with the automorphism groups of various 
Lie algebras of vector fields on $\AA^n$. 
Below $\VF^0(\AA^n)$ 
stands for the Lie algebra of vector fields 
with zero divergence. 
\begin{thm}
There are  isomorphisms 
\[\Aut(\AA^n)\cong\Aut_{\rm Lie}(\VF(\AA^n))\cong 
\Aut_{\rm Lie}(\VF^{\mathrm c}(\AA^n))
\cong\Aut_{\rm Lie}(\VF^0(\AA^n))\]
induced by the natural action of the group $\Aut(\AA^n)$ 
on the Lie algebra of vector fields $\VF(\AA^n)$. 
\end{thm}
The first isomorphism was obtained
in \cite[Corollary]{Rud86}, \cite[Sec. 4]{Kul93} 
and  \cite{Bav17a} 
along different approaches;
in \cite{Kul93} it appeared in relation 
with the Jacobian Conjecture. 
The other isomorphisms were established
independently
in \cite{Bav17b} 
and \cite{KrRe17}.  See \cite{Bav14} for 
the automorphism group of the Lie algebra 
of triangular derivations. 

See also \cite[Th\'eor\`eme A]{CL06} 
for a description of the lattices in simple real Lie groups 
that can be realized as subgroups of
$\Aut(\AA^2_{\CC})$. Furthermore, see
\cite{AGL24} 
and the references therein 
for studies on (finite dimensional)  
Lie algebras of derivations of polynomial rings. 
\subsection{Tamizable automorphisms}
\begin{defn}[{\rm \cite[Definition 4.1]{MP09}}]
 An automorphism $\phi\in\Aut(\AA^n)$ is 
 called \emph{tamizable}
(resp., \emph{linearizable, triangulable}) 
if it is conjugate in 
$\Aut(\AA^n)$ to a tame 
(resp. linear, triangular) automorphism.
\end{defn}
\begin{que}[{\rm \cite{MP09}}]
Is it true that any $\phi\in\Aut(\AA^n)$ is tamizable? 
In particular, is the Nagata automorphism 
$\Phi$ tamizable?  
\end{que}
\begin{exa}
It is known that the Nagata automorphism $\Phi$ 
is not conjugate to a triangular one, 
see \cite{Bas84}.
Nevertheless, it is \emph{shifted linearizable}. 
More precisely, 
$2{\rm id}\cdot\Phi$ is conjugate to $2{\rm id}$ 
in $\Aut(\AA^3)$, 
see \cite[Theorem 3.3]{MP09}; 
cf. also \cite[Lemma 15.11.1]{FuKr18} for a general 
result of this type.
\end{exa}
As follows from the Jung-van der Kulk Theorem 
\ref{thm:JvdK}, $\Aut(\AA^2)$ 
is algebraically generated, see subsection \ref{ss:ags} 
for the definition.
It is known that $\Aut(\AA^n)$ is connected for any $n$, 
see  \cite[Lemma 4]{Sha81}. 
The following natural question arises: 
\begin{que} Is $\Aut(\AA^n)$ algebraically generated 
for $n\ge 3$?
\end{que}
The locally finite automorphisms of $\Aut(\AA^n)$
were studied in \cite{FuMa07}. It is known that the subgroup 
$\langle{\rm LF}_n\rangle$ of $\Aut(\AA^n)$ 
generated by all locally finite automorphisms is normal, 
see \cite{Zyg11}. By  
the Jung-van der Kulk Theorem \ref{thm:JvdK} we have
$\langle{\rm LF}_2\rangle=\Aut(\AA^2)$.  For $n=3$, 
$\langle{\rm LF}_3\rangle$
contains the Nagata automorphism.
However, it is unknown whether 
$\Aut(\AA^3)=\langle{\rm LF}_3\rangle$
(\cite{Zyg11}; cf. also subsection \ref{ss:cotame}).
\subsection{Stable tameness} 
\begin{defn} 
A wild automorphism $\phi\in\Aut(\AA^n)$ is called 
\emph{stably tame} if 
$\phi\times {\rm \id}_{\AA^\kappa}\in\Aut(\AA^{n+\kappa})$ 
is tame  for some 
$\kappa\ge 1$. 
The  minimal number $\kappa$ 
of additional variables 
will be called 
\emph{excessive dimension}.
\end{defn}
\begin{thm}[{\rm \cite{Smi89}}] 
The Nagata automorphism $\Phi$ 
is stably tame with excessive dimension 
$\kappa=1$. 
\end{thm}
By \cite{Spo07} the same is true for the 
Anick automorphism of $\AA^4$, 
cf. Remark \ref{rem:Anick}.3.
More generally, the following holds. 
\begin{thm}[{\rm \cite[Corollary 4.9]{BvdEW12}}]
Let as before $\G\subset\Aut(\AA^3)$ 
be the subgroup of automorphisms fixing $z$.
Then any element $\phi\in\G$ is stably tame 
with excessive dimension $\kappa\le 3$.
\end{thm}
A similar result is also known for $\Aut(A_3)$ 
with $\kappa=1$, see \cite[Theorem 1.1]{BKY12b}. 
\subsection{Cotame automorphisms}\label{ss:cotame}
\begin{defn}[\emph{Cotame automorphisms}; 
cf.\ \cite{EL19}]\label{sit:tame} 
One says that $h\in\Aut(\AA^n)\setminus  {\Aff}_n$ 
is \emph{cotame} 
if $\langle {\Aff}_n,h\rangle\supset {\Tame}_n$ 
and \emph{topologically cotame} if 
$\overline{\langle {\Aff}_n,h\rangle}
\supset {\Tame}_n$. 
\end{defn}
\begin{thm}[{\rm \cite{Bod02}, \cite[1.4 and 
Theorem 1.8]{Bod05}, \cite{EL19}}]
\label{rem:cotame} 
For $n\ge 3$  
any triangular  $h\in \Jonq_n\setminus\Aff_n$ 
is cotame, 
while there is no triangular cotame 
$h\in \Aut(\AA^2)\setminus\Aff_2$.
\end{thm}
The following theorem of Edo
extends and refines the earlier results of  
Bodnarchuk (\cite[Thm.\ 3]{Bod01}) and Furter 
(\cite[Thm.\ D]{Fur15}); cf. also 
\cite[Theorem 1.4]{FuPo18}.
\begin{thm}[{\rm \cite[Thm.\ 1.2]{Edo18}}]
\label{thm:I} 
For $n\ge 2$ any element $h\in 
\Aut(\AA^n)\setminus\Aff_n$ 
is topologically cotame. 
\end{thm}
\subsection{Amalgamation}
By Jung-van der Kulk's Theorem \ref{thm:JvdK}, 
$\Aut(\AA^2)$ 
is an amalgam of two subgroups along their 
intersection.
A similar decomposition exists for the group 
$\GL(2, \kk[x])$,
see \cite{Nag59}, and for the group $\Aut^0(A_2)$ 
of augmentation-preserving automorphisms 
of the free associative algebra $A_2$, see 
\cite[Theorem 1.1]{SY98}.

It is unknown, however, whether a similar result
holds for  $\Aut(\AA^n)$ 
when $n\ge 3$. For the tame subgroup 
$\Tame_3\subset\Aut(\AA_3)$ we have 
the following result.
\begin{thm}[{\rm\cite{Wri15}, see also 
\cite{Wri17}}]
$\Tame_3$  can be realized as a generalized 
amalgamated product of three subgroups
amalgamated along pairwise intersections. 
\end{thm}
See \cite[Theorem 1]{ANU20} and the 
references therein 
for analogues of the Jung-van der Kulk's Theorem
for  the automorphism group of the free 
Lie algebra in three variables and 
the subgroup of tame automorphisms. 
Similar results are also known to hold for 
 certain $\Aut$-homogeneous affine surfaces 
among the so-called \emph{Danilov-Gizatullin surfaces}. 
\begin{exa}[\emph{Danilov-Gizatullin surfaces}] 
\label{ex:DG}
These are smooth affine surfaces 
$X_d= F_n \setminus S_d$ 
where $\pi_n\colon F_n\to\PP^1$ is the $n$th 
(smooth projective) Hirzebruch ruled surface and 
the curve $S_d\subset F_n$ 
is an ample section of $\pi_n$ with $S_d^2=d> n$.
Two such affine surfaces are isomorphic if and
only if they share the same invariant $d$,  
see \cite[Theorem 5.8.1]{GD77}. 
The group $\Aut(X_d)$
 has exactly $\lfloor d/2\rfloor$ 
 conjugacy classes of
maximal tori of rank 1, 
see \cite[Prop. 5.15]{FKZ07}. 
The group $\Aut(X_d)$ is an amalgam of 2 
nested ind-subgroups if $d = 3$ 
and of 3 nested ind-subgroups if $d = 4,5$,
see \cite[Sections 5-8]{GD77}. 
Assume that the base 
field $\kk$ is uncountable. Then
$\Aut(X_6)$ is not an amalgam of 
two nested ind-groups. 
For $d\ge 7$,
$\Aut(X_d)$ is not an amalgam of 
any countable number of nested ind-groups, see 
\cite[Corollary 5.6]{KPZ17}.

One can find similar results for other 
$\Aut$-quasihomogeneous affine surfaces 
in \cite{GD75, GD77}, see also \cite[Sect. 5]{KPZ17}. 
\end{exa}
\subsection{Linearization of 
reductive group actions} 
Let $G$ be a reductive linear algebraic group 
acting effectively on $\AA^n$
with a fixed point $x$. Then the 
tangent representation 
of $G_x$ on $T_xX$ is faithful.
Luna's \'Etale Slice Theorem 
\cite{Lun73} guarantees 
a local linearization of the $G_x$-action 
near $x$. 
This justifies an expectation that any action of 
a reductive group on $\AA^n$
can be linearized, that is, is conjugate in 
$\Aut(\AA^n)$ to a linear action; 
see e.g. \cite[Conjecture 3.1]{Kam79} 
and \cite[Section 9.4]{vdE00}. 
However, the latter expectation turns out 
to be wrong as we will see below in \ref{counterexamples}.

\par\medskip
The first results in this context are due  
to Bia\l{}ynicki-Birula who proved the following 
in \cite{Bia66-67}.
\begin{thm}
A faithful action of a torus of dimension 
$\geq n-1$ on $\AA^n$ is linearizable.
In particular, all maximal tori in $\Aut(\AA^n)$ 
are pairwise conjugate. 
\end{thm}
The latter fails for affine surfaces and 
affine threefolds, in general,
see Examples \ref{ex:Dani} and \ref{ex:DG}.

It then took a long time to settle the next open problem, 
namely the linearizability of  a
$\mathbb{G}_{\mathrm m}$-action on $\AA^3_{\CC}$. 
This was done in a series of fundamental papers 
of Koras and Russell,
see e. g. 
\cite{KoRu97}, 
completed finally in \cite{KKMLR97}.
\begin{thm}
Any action of $\Gm$ on $\AA^3_\CC$  
is linearizable.
\end{thm}
The proof uses a classification of certain
Koras-Russell affine  threefolds diffeomorphic to 
$\mathbb{R}^6$ and non-isomorphic 
to $\AA^3_\CC$, see  \cite{KML97}.

Certain actions of codimension 
2 tori can be linearized also in higher dimensions, 
see  \cite{KamRu82} for (a) 
and \cite{KoRu89} for (b). 
\begin{thm} Assume that 
the algebraic torus $T= (\Gm)^{n-2}$
acts effectively on $\AA^n$. 
Then this action can be linearized
provided one of the following holds:
\begin{enumerate}
\item[(a)]
The action is "unmixed," 
 that is, no nontrivial character, 
 together with its inverse, 
 occur as weights of semi-invariant polynomials;
\item[(b)]
the fixed point set of $T$ has positive dimension. 
\end{enumerate}
\end{thm}

Next we look at the case of the affine plane  $\AA^2$.
Using the amalgamated product structure of 
$\Aut(\AA^2)$ and a famous result of Serre's 
\cite[I.4.3.8]{Ser03} on subgroups of 
bounded length in an amalgamated product,
one gets the following theorem. 
\begin{thm}[{\rm \cite[Theorem~4.3]{Kam79}}]
\label{Kambayashi.thm}
An algebraic subgroup $G$ of $\Aut(\AA^2)$
is conjugate to a subgroup of either  
$\Aff_2$ or $\Jonq_2$. In particular, 
every reductive group action on $\AA^2$ 
is linearizable.
\end{thm}
In particular, 
any action of a finite group on 
$\AA^2$ is conjugate to a linear action. 
This was first proven by geometric means by
T. Igarashi (1975, unpublished), 
see also \cite[Theorem~2]{Fur83} and
\cite[Theorem 2.4.6]{Miy24}.
\par\medskip
For actions of reductive groups 
we have the following consequence of Luna's 
\'Etale Slice Theorem, see \cite[III. Corollaire~2]{Lun73}; 
cf. also \cite[Corollary 10.7]{BH85}.
\begin{thm}
Let a reductive algebraic group $G$ act 
on an affine variety $X$. 
Assume that every  $G$-invariant function is 
constant and that $G$ has 
a smooth fixed point in $X$.
Then $X$ is $G$-equivariantly 
isomorphic to a representation of $G$.
\end{thm}
Applied to the actions on $\AA^n$, 
this gives the following corollary.
\begin{cor}[\rm{\cite[Proposition 5.1]{KrPo85}}]
If a reductive group acts on $\AA^n$ 
with a dense orbit, then the action is linearizable.
\end{cor}
For semisimple groups we have the following results, 
see  \cite{KrPo85} and \cite{KrRu14} for $n=3$ 
and \cite{Pan84-86} for $n=4$.
\begin{thm}
For $n\le 4$, any action of a semisimple algebraic group 
on $\AA^n$ is linearizable. 
\end{thm}
In the paper \cite{KrSc92} Kraft and Schwarz study 
the linearization problem for reductive group actions 
on smooth acyclic affine varieties $X$ assuming that 
the quotient $X/\!/ G$ has dimension 1. It turns out that 
there is always a fixed point $x_0 \in X$ 
(\cite[Theorem 1]{KrSc92}) 
so that the tangent representation 
$V := T_{x_0} X$ is the candidate 
for the linearization. 
One of the main results is the following 
(\cite[Theorem 5]{KrSc92}).
\begin{thm} $X$ is $G$-isomorphic to $V:=T_{x_0}X$ in the following cases:
\be
\item $V$ is a semifree $G$-module;
\item $G$ is a torus;
\item $\dim  V^{G^0} = 1$;
\item $\dim X \leq 3$;
\item $G^0$ is semisimple;
\item $V$ is selfdual as a $G^0$-module.
\ee
\end{thm}
As a corollary one gets the next result, cf. \cite[Theorem~A]{KrRu14}.
\begin{cor}
Every faithful action of a non-finite reductive group 
on $\AA^3$ is linearizable.
\end{cor}
As for linearizable finite group actions, 
Kraft and Procesi suggested 
a conjectural linearization procedure that uses 
the Reynolds operator.  
This  linearization procedure works indeed
in the case of involutions of $\AA^n$ 
acting with a fixed point, provided the 
associated Taylor decomposition 
of the Reynolds average is of bounded length, 
see \cite{Jur90}.
\subsubsection*{Counterexamples}
\label{counterexamples}
The first examples of non-linearizable actions 
of connected 
reductive groups on $\AA^n$ are due to 
Schwarz \cite{Sch89}.
These are 
an ${\rm O}_2(\CC)$-action on $\AA_{\CC}^4$ 
and an $\SL_2(\CC)$-action on $\AA_{\CC}^7$. 
In a sense, the non-linearizability is 
a general phenomenon. 
\begin{thm}[{\rm \cite{Kno91}}]
For any  non-commutative connected reductive 
algebraic group $G$
 there exists a non-linearizable $G$-action 
 on some affine space $\AA^n$. 
\end{thm}
Many more examples are given in 
\cite[VII.5. Examples]{KrSc92}.
\par\medskip
All known examples of non-linearizable actions 
on $\AA^n$ are obtained from 
non-trivial algebraic $G$-vector bundles 
over $G$-modules,
following an idea from 
\cite{BH87}, see also 
\cite{Kra89}, \cite{KrSc92}, \cite{MP91} 
and \cite{MMP91}. 
Indeed, by the Quillen-Suslin theorem, 
such a bundle is trivial when  ignoring 
the $G$-action.
Hence, its total space is isomorphic to 
some $\AA^n$. 
The induced $G$-action on $\AA^n$ 
is non-linearizable due to the fact that
our bundle is non-trivial in the category
of $G$-vector bundles. 

However, this method does not work for 
an abelian reductive group $G$.
Indeed, in this case every $G$-vector 
bundle over 
a $G$-module is trivial, see \cite{MMP96}.
The following questions 
remain open, 
see e.g. \cite{KrSc95} and \cite{Mas03}:
\begin{que} Is there a non-linearizable action 
of an abelian group on 
$\AA^n$, for instance, of a cyclic group 
$\ZZ/p\ZZ$ on $\AA^n$, $n\ge 3$ and of 
$\Gm$ on $\AA^n$, $n\ge 4$?
\end{que}
Notice that 
there is an example of Asanuma of a 
non-linearizable action of 
$\Gm(\mathbb{R})$ on $\AA^n_{\mathbb{R}}$,
see \cite{Asa99}. 

\par\medskip
There are several examples of 
non-linearizable actions 
of non-abelian finite groups on $\AA^n$, for instance, 
of the dihedral groups $D_{10}$, $D_{14}$ 
and of of the symmetric group 
$S_3$ on $\AA^4$, 
see \cite[Theorem 5.1]{MJ94} and  \cite{FM15}; 
cf. also \cite{MP95}, \cite{MMP96}, \cite{Mas03}.
It seems that no example of this kind is known for 
a finite group action on  $\AA^3$, see e.g. 
\cite{FM15}. See also \cite{PeRa86} for the  
Lefschetz fixed-point formula in the context 
of finite group actions 
on affine varieties.

The linearization problem is closely related to the 
Abhyankar--Sathaye problem on linearization of 
embeddings $\AA^k\hookrightarrow\AA^n$ 
and on variables in polynomial rings. 
We send the reader e.g. to \cite{DG15},
 \cite{Kal09}, \cite{KrRu14}, 
\cite{vdEMV07} 
and the literature therein for further reading. 
\subsection{Solvable and unipotent ind-subgroups: 
triangulation}\label{ss:triangulation}
We start with the following result.
\begin{prop}[{\rm \cite[Proposition 10.5.1]{FuKr18}}] 
Let $U \subset \Aut(X)$ be
a commutative unipotent algebraic subgroup 
of dimension $n$. 
Assume that ${\rm Cent}_{\Aut(X)}(U)=U$.
Then for any $x \in X$ the orbit map 
$U \to X$, $u \mapsto ux$, 
is an isomorphism. In particular, $X$ 
is isomorphic to $\AA^n$.
\end{prop}
The unipotent elements of the de Jonqui\`eres group 
$\Jonq_n$ form 
a solvable subgroup $\JJJ_n$ of derived length 
$n$,  see \cite[Lemma 3.2]{FuPo18}. 
This subgroup is not nilpotent for $n>1$; 
see Remark \ref{rem2} above.
\begin{defn} We say that a subgroup 
$G\subset\Aut(\AA^n)$ is \emph{triangulable} 
if it is conjugate in $\Aut(\AA^n)$ 
to a subgroup of $\Jonq_n$. It is said to be 
\emph{stably triangulable} if  $G\times \id_{\AA^k}$ 
is triangulable in 
$\Aut(\AA^{n+k})$ for some $k\ge 0$. 
\end{defn}
By \cite{Ren68}
any unipotent subgroup of 
$\Aut(\AA^2)$ is triangulable, 
cf. Theorem~\ref{Kambayashi.thm}. 
This is not any longer  true in higher dimensions.
However, the next theorem provides
 a large class of triangulable unipotent subgroups 
 in any dimension.
\begin{thm}[{\rm\cite[Theorem D]{KrZa24}}]
\label{thm:triang}
Let $U \subset \Aut(\AA^n)$ be a nested 
unipotent subgroup.  
If $U$ acts transitively on $\AA^n$, 
then $U$ is triangulable.
\end{thm}
A similar result holds for Lie subalgebras 
of vector fields.
\begin{thm}[{\rm\cite[Theorem 6]{Sku21}}] 
Let $L\subset \Lie(\Aut(\AA^n))$ 
be a Lie subalgebra that consists 
of locally nilpotent derivations. 
Assume that 
$\bigcap_{\partial\in L}\ker(\partial)=\kk$. 
Then $L$ is ${\rm Ad}$-conjugate 
to a subalgebra of 
the Lie algebra $\Lie(\Jonq_n)$ of 
triangular derivations. 
\end{thm}
The transitivity assumption in Theorem 
\ref{thm:triang} is essential, 
as the following example shows. 
\begin{exa}\label{ex:NB} 
Consider the Nagata automorphism 
$\Phi=\exp(\sigma\p)$, see subsection \ref{ss:tame}.
It is contained in the $\Ga$-subgroup 
$U=\exp(t\sigma\p)$ of $\Aut(\AA^3)$.
According to Bass \cite{Bas84}, $U$ is
 not triangulable. 
Indeed, the fixed points set $(\AA^3)^{\Ga}$ 
is a hypersurface 
with an isolated singularity. 
However,  the fixed point set of a triangulable 
$\Ga$-subgroup
is isomorphic to a cylinder $X \times\AA^1$ 
which cannot have an isolated singularity. 
Using the same idea one can construct  
non-triangulable 
$\Ga$-subgroups of any $\Aut(\AA^n)$, $n\ge 3$,
see \cite{Pop87}. 
\end{exa}
According to Freudenburg 
\cite[Lemma 3.36]{Fre17}, 
the Nagata-Bass $\Ga$-subgroup 
$U$ is not stably triangulable in $\Aut(\AA^4)$. 
In fact, no $\Ga$-subgroup of Bass-Popov type 
is stably triangulable,
whatever is the number of additional variables.   
\begin{prop}[{\rm\cite[Proposition 1.4.1]{KrZa24}}]
\label{Bass.prop}
Consider a $\Ga$-sub\-group 
$U\subset\Aut(\AA^n)$. 
Assume that the fixed point set $(\AA^n)^U$ 
is a hypersurface 
with an isolated singularity. 
Then $U$ is not stably triangulable.
\end{prop}
There is, however, the following positive result. 
\begin{thm}[{\rm\cite{Kal04}}]
Assume that $\kk=\CC$, and let a 
$\Ga$-subgroup $U\subset\Aut(\AA^3)$
acts freely on $\AA^3$. Then $U$ is conjugate 
in $\Aut(\AA^3)$ 
to a subgroup of translations. In particular, 
$U$ is triangulable. 
\end{thm}
An analogue of the former assertion fails in $\AA^4$, 
see \cite{Win90}.
For further results in this direction, see e.g. 
the survey article \cite{Kal09} and
 \cite{GMM21} for $\Ga$-actions on affine varieties,
\cite{Fre17} and \cite{vdEKC21} for $\Ga$-actions 
on affine spaces and Hilbert's 14th problem. See also
\cite{Pop17} and the references therein for (stably) 
birationally triangulable and non-triangulable actions
of unipotent groups. 
\section{Tits' type alternative for automorphism groups}
A group $G$ is called \emph{virtually solvable} 
resp.,  \emph{virtually nilpotent}, 
etc. if $G$ is a finite extension of a solvable resp. 
nilpotent, etc. group. 
The classical {\it Tits alternative\/} says that any linear group 
$G$, that is 
a subgroup of $\GL(n,\kk)$, is either virtually solvable, 
or contains 
a free group $\FF_2$ on two generators. 
For finitely generated subgroups, 
the Tits alternative holds 
without any restriction on the base field.

There are several similar statements 
for automorphism groups. For instance, an analogue of
 the Tits alternative holds for $\Aut(\AA^2_{\CC})$, 
see  \cite[Theorem 2.4]{Lam01} or 
\cite[Proposition 20.17]{Lam24}.

Notice that $\Aut(\AA^2_{\kk})$ 
 for an infinite field $\kk$
is not isomorphic to a linear group 
over any extension of $\kk$, and the same 
holds for its de Jonqi\`eres subgroup $\Jonq_2$, 
see \cite[Proposition 2.3 and Remark 2.4]{Cor17} 
and also \cite{Mat23}. 
However,  by \cite[Theorem A]{Mat23}
the codimension 6 normal subgroup 
\[G=\{\phi\in \Aut(\AA^2_{\kk})\,|
\,\phi(0)=0,\,\,\,d\phi(0)=\id\}\]
of $\Aut(\AA^2_{\kk})$ 
is a linear group over a suitable  field extension.

The unipotent triangular 
subgroup $\JJJ_n$ 
of $\Aut(\AA^n_{\kk})$ over an arbitrary field $\kk$ 
is a semidirect product of abelian groups, 
see \cite[Theorem~1]{BNS12}.
For $n\ge 3$ this group
is nonlinear (it cannot be represented 
by matrices over any field), 
see \cite[Theorem~2]{Sos07}. 
The same holds for the tame subgroup 
$\Tame_n=\Tame(\Aut(\AA^n_{\kk}))$ provided 
${\rm char}(\kk)=0$ and $n\ge 4$, 
see \cite{RCS04}. Similar results are also true for 
the automorphism groups of free associative 
algebras, see {\it ibid}.

For the proof of the following theorem 
see \cite{Can11} 
and \cite{Ure21} for (a), 
\cite[Theorem C]{BFL14} for (b) 
and \cite[Main Theorem]{LP22} for (c). 
 In Theorems \ref{thm:Tits1} 
and \ref{thm:Tits2} 
we assume that 
the base field $\kk$ is 
algebraically closed of characteristic zero.
\begin{thm}\label{thm:Tits1}  The Tits 
alternative holds for the following groups:
\begin{itemize}
\item[(a)] the birational automorphism group of any 
compact complex K\"{a}hler surface;
\item[(b)] $\Tame(\Aut(\SL(2,\kk)))$ 
(see Remark \ref{rem:Anick}.3);
\item[(c)] $\Tame_3=
\Tame(\Aut(\AA^3_{\kk}))$.
\end{itemize}
\end{thm}
In the case of algebraically generated 
groups of automorphisms, 
we have the following Tits' type alternatives.
\begin{thm}[{\rm\cite[Theorem 1.1]{AZ22} and 
\cite[Theorem 6.1]{AZ23}}]\label{thm:Tits2} $\,$
\begin{itemize}
\item[(a)] Let $X$ be an affine algebraic surface over $\kk$
and $G\subset\Aut(X)$ be a subgroup generated 
by $\Ga$-subgroups $U_1,\ldots,U_k$. 
Then either $G$ contains $\FF_2$ 
or $G$ is a metabelian unipotent algebraic group.
\item[(b)] Let $X$ be a toric affine variety  over $\kk$
and $G\subset\Aut(X)$ 
be a subgroup
generated by 
Demazure root subgroups $U_1,\ldots,U_k$,
see Definition \ref{Demazure}.
Then either $G$ contains $\FF_2$ or $G$ is a 
unipotent algebraic group.
\end{itemize}
\end{thm}
Similar results also hold  for subgroups of 
$\Tame(\Aut(\SL(2,\kk)))$
and $\Tame(\Aut(\AA^3_{\kk}))$, 
see \cite[Corollary 6.4]{AZ23}.

\medskip

\noindent {\bf Acknowledgments.}  
We would like to thank 
Ivan Arzhantsev, 
Alexander Perepechko and
Andriy Regeta for their comments and references. 
 Our thanks also go to the anonymous referee
for useful remarks. 

\renewcommand{\MR}[1]{}


\bigskip
\begin{thebibliography}{KKMLR97}
%
\bibitem[ANU20]{ANU20}
Alibek~A.~Alimbaev, Ruslan~Zh.~Nauryzbaev, 
and Ualbai~U.~Umirbaev,
\emph{On the automorphisms 
of a free Lie algebra 
of rank $3$ over an integral domain},
Sib. Math. J. 61:1 (2020), 1--10. 
%
\bibitem[And19]{And19} 
Rafael~B.~Andrist, 
\emph{Integrable generators of Lie algebras 
of vector fields on $\CC^n$},
Forum~Math.\ 31 (2019), 943--950.
%
\bibitem[Arz18a]{Arz18a}
Ivan~Arzhantsev,
\emph{Infinite transitivity and special automorphisms}, 
Ark. Mat. 56:1 (2018), 1--14. 
%
\bibitem[Arz18b]{Arz18b}
Ivan~Arzhantsev,
\emph{Gale duality and homogeneous toric varieties}, 
Comm. Algebra 46:8 (2018), 3539--3552.
%
\bibitem[Arz23]{Arz23}
Ivan~Arzhantsev,
\emph{Automorphisms of algebraic varieties 
and infinite transitivity},
St. Petersburg Math. J. 34:2 (2023), 143--178.
%
\bibitem[AFK$^+$13a]{AFK$^+$13a}
Ivan~Arzhantsev, Hubert~Flenner, 
Shulim~Kaliman, Frank Kutzschebauch, and
Mikhail~Zaidenberg, 
\emph{Flexible varieties and automorphism groups}, 
Duke Math. J. 162:4 (2013), 767--823.
%
\bibitem[AFK$^+$13b]{AFK$^+$13b} 
 Ivan~Arzhantsev, Hubert~Flenner, 
Shulim~Kaliman, Frank Kutzschebauch, and
Mikhail~Zaidenberg, 
\emph{Infinite transitivity on affine varieties}, 
In:
Birational Geometry, Rational Curves, 
and Arithmetic. Fedor~Bogomolov, 
Brendan~Hassett, Yuri~Tschinkel (eds.), 1--14.
Springer-Verlag, New York e.a. 2013.
%
\bibitem[AG10]{AG10}
Ivan~V.~Arzhantsev  and Sergei~A.~Gaifullin, 
\emph{Cox rings, semigroups and automorphisms of 
affine algebraic varieties}, 
Sb. Math.  201:1 (2010), 1--21.
%
\bibitem[AG17]{AG17}
Ivan~Arzhantsev and Sergey~Gaifullin,
\emph{The automorphism group of a rigid affine variety},
Math. Nachr. 290:5–6  (2017), 662--671. 
%
\bibitem[AGL24]{AGL24}
Ivan~Arzhantsev Sergey~Gaifullin,  
and Viktor~Lopatkin,
\emph{On finite-dimensional homogeneous 
Lie algebras of derivations of polynomial rings},
arXiv:2408.05627 (2024).
%
\bibitem[AKZ12]{AKZ12}
Ivan~Arzhantsev, Karine~Kuyumzhiyan, 
and Mikhail~Zaidenberg, 
\emph{Flag varieties, toric varieties, and suspensions: 
Three instances of infinite transitivity}, 
Sb. Math. 203 (2012), 923--949. 
%
\bibitem[AKZ19]{AKZ19}
Ivan~Arzhantsev, Karine~Kuyumzhiyan, 
and Mikhail~Zaidenberg, 
\emph{Infinite transitivity, finite generation, 
and Demazure roots},
Adv. Math. 351 (2019), 1--32. 
%
\bibitem[APS14]{APS14}
Ivan~Arzhantsev, Alexander~Perepechko, 
and Hendrik~S\"{u}ss, 
\emph{Infinite transitivity on universal torsors}, 
J. London Math. Soc. (2) 89:3 (2014), 762--778.
%
\bibitem[AZ22]{AZ22} 
Ivan~Arzhantsev and Mikhail~Zaidenberg, 
Tits-type alternative for groups acting 
on toric affine varieties, 
Int. Math. Res. Not. 11  (2022), 8162--8195.
%
\bibitem[AZ23]{AZ23}
Ivan~Arzhantsev and Mikhail~Zaidenberg, 
\emph{Tits-type alternative for certain groups acting 
on algebraic surfaces},
Proc. Amer. Math. Soc. 151:7 (2023),  2813--2829.
%
\bibitem[ArZa23]{ArZa23}
Ivan~Arzhantsev and Yulia~Zaitseva,
\emph{Affine homogeneous varieties 
and suspensions}, Res. Math. Sci. 11:2 (2024), 
article 27.
%
\bibitem[Asa99]{Asa99}
Teruo~Asanuma, 
\emph{Non-linearizable algebraic 
$k^*$-actions on affine spaces},
Invent. Math., 138:2 (1999), 281--306.
%
\bibitem[Bas84]{Bas84}
Hyman~Bass, 
\emph{A non-triangular action of 
$\mathbb{G}_{\mathrm a}$ on $\AA^3$}, 
J. Pure Appl. Algebra 33 (1984), 1--5.
%
\bibitem[BH85]{BH85}
Hyman~Bass and William~J.~Haboush, 
\emph{Linearizing certain reductive 
group actions}, 
Trans. Amer. Math. Soc. 292 (1985), 463--482.
%
\bibitem[BH87]{BH87}
Hyman~Bass and William~J.~Haboush, 
\emph{Some equivariant 
{$K$}-theory of affine algebraic
group actions}, 
Comm. Algebra \textbf{15} (1987), no.~1-2, 181--217.
%
\bibitem[BL83]{BL83}
Hyman~Bass and Alexander~Lubotzky, 
\emph{Automorphisms of groups 
and of schemes of finite type}, 
Israel J. Math. 44 (1983), 1--22.
%
\bibitem[BNS12]{BNS12}
Valeriy~G.~Bardakov, Mikhail~V.~Neshchadim, 
and Yury~V.~Sosnovsky,
\emph{Groups of triangular automorphisms 
of a free associative algebra 
and a polynomial algebra}, 
J. Algebra 362 (2012), 201--220.
%
\bibitem[Bav14]{Bav14}
Vladimir V. Bavula, 
\emph{The groups of automorphisms of 
the Lie algebras of triangular polynomial derivations}, 
J. Pure Appl. Algebra 218 (2014), 829--851.
%
\bibitem[Bav17a]{Bav17a}
Vladimir~V.~Bavula, 
\emph{The group of automorphisms of 
the Lie algebra of derivations of a polynomial algebra},
 J. Algebra Appl. 16:5 (2017), 1750088. 
%
 \bibitem[Bav17b]{Bav17b}
Vladimir~V.~Bavula, 
\emph{The groups of automorphisms of 
the Lie algebras of polynomial vector fields 
with zero or constant divergence}, 
Commun. Algebra 45:3 (2017), 1114--1133.
%
 \bibitem[BKYE16]{BKYE16}
Alexei~Belov-Kanel, Jie-Tai~Yu, and  Andrey Elishev,
\emph{On the Zariski topology of automorphism 
groups of affine spaces and algebras},
arXiv:1207.2045v6 (2016).
%
\bibitem[BKY12b]{BKY12b}
Alexei~Belov-Kanel and Jie-Tai~Yu,
\emph{Stable tameness of automorphisms of 
$F\langle x, y, z\rangle$ fixing $z$}, 
Sel. Math. New Ser. 18 (2012), 799--802. 
%
\bibitem[BH03]{BH03}
Florian~Berchtold and J\"{u}rgen~Hausen, 
\emph{Demushkin's Theorem in 
codimension one}, 
Math. Z. 244 (2003),  697--703. 
%
\bibitem[BW00]{BW00}
Yuri~Berest and George~Wilson, 
\emph{Automorphisms and ideals of the Weyl algebra}, 
Math. Ann. 318:1 (2000), 127--147. 
%
\bibitem[BvdEW12]{BvdEW12}
Joost~Berson, Arno van den Essen, 
and David Wright,
\emph{Stable tameness of two-dimensional 
polynomial automorphisms over a regular ring},
Adv. Math. 230 (2012), 2176--2197.
%
\bibitem[BPZ24]{BPZ24}
Oksana~Bezushchak, Anatoliy~Petravchuk, 
and Efim~Zelmanov,
\emph{Automorphisms and derivations of 
affine commutative and PI-algebras},
Trans. Amer. Math. Soc. 377:2  (2024), 
1335--1356.
%
\bibitem[Bia66-67]{Bia66-67}
Andrzej~Bia\l{}ynicki-Birula, 
\emph{Remarks on the action of 
an algebraic torus on $k^n$. I, II}, 
Bull. Acad. Polon. Sci. Ser. Sci. Math. 
Astronom. Phys. 
14 (1966), 177--181; ibid. 15 (1967), 123--125.
%
\bibitem[BFL14]{BFL14}
Cinzia~Bisi, Jean-Philippe~Furter, 
and St\'ephane~Lamy,
\emph{The tame automorphism group 
of an affine quadric threefold acting 
on a square complex},
J. Ec. polytech. Math. 1 (2014), 161--223. 
%
\bibitem[BD15]{BD15}
J\'er\'emy~Blanc and Adrien~Dubouloz,
\emph{Affine surfaces with a 
huge group of automorphisms},
Int. Math. Res. Not. 2 (2015), 422--459. 
%
\bibitem[BG23]{BG23}
Viktoria~Borovik and Sergey~Gaifullin,
\emph{Isolated torus invariants and 
automorphism groups of rigid varieties},
arXiv:2007.07882 (2023).
%
\bibitem[Bod01]{Bod01} Yuri~Bodnarchuk, 
\emph{Some extreme properties 
of the affine group as 
an automorphisms group of the affine space}, 
Contribution to General Algebra 13 (2001), 15--29.
%
\bibitem[Bod02]{Bod02} Yuri~Bodnarchuk, 
\emph{Affine group as a subgroup of 
biregular transformation
group of an affine space},  (Ukrainian). 
NaUKMA Zap., Ky\"{i}v 20 (2002), 6--10.
%
\bibitem[Bod05]{Bod05}
Yuri~Bodnarchuk, 
\emph{On generators of the tame invertible
polynomial maps group},
Int.\ J.\ Algebra Comput.\ 15 (2005), 851--867.
%
\bibitem[Can11]{Can11}
Serge~Cantat, 
\emph{Sur les groupes de transformations 
birationnelles des surfaces}, 
Ann. Math. (2) 174:1 (2011), 299--340. 
%
\bibitem[BG23]{BG23}
Viktoriia Borovik and Sergey Gaifullin,
\emph{Isolated torus invariants and 
automorphism groups of rigid varieties,
 arXiv:2007.07882v2 (2023).
%
\bibitem[CL06]{CL06}
Serge~Cantat and St\'ephane Lamy,
\emph{Groupes d'automorphismes polynomiaux du plan},
Geom. Dedicata 123  (2006), 201--221. 
}
%
\bibitem[CRX19]{CRX19}
Serge~Cantat, Andriy~Regeta and Junyi~Xie,
\emph{Families of commuting automorphisms, 
and a characterization of the affine space},
Amer. J. Math. 145:2 (2019), 413--434.
%
\bibitem[CPPZ21]{CPPZ21}
Ivan~Cheltsov, Jihun~Park, Yuri~Prokhorov, 
and Mikhail~Zaidenberg,
\emph{Cylinders in Fano varieties},
EMS Surv. Math. Sci. 8 (2021), 39--105,
%
\bibitem[ChDu23]{ChDu23}
Mikhael~Chitayat and Adrian~Dubouloz,
\emph{The rigid Pham-Brieskorn threefolds},
arXiv:2312.07587 (2023).
%
\bibitem[CD03]{CD03}
Arjeh~M.~Cohen and Jan~Draisma, 
\emph{From {L}ie algebras of vector fields to
algebraic group actions}, Transform. Groups  
8:1 (2003), 51--68.
%
\bibitem[Coh85]{Coh85}
Paul~M.~Cohn, 
\emph{Free rings and their relations}, 
2nd ed., Academic Press, London 1985.
%
\bibitem[Cor17]{Cor17}
 Yves~Cornulier, 
 \emph{Nonlinearity of some subgroups 
 of the planar Cremona group}, 
 arXiv:1701.00275 (2017).
%
\bibitem[Cze71-72]{Cze71-72}
Anastasia~J.~Czerniakiewicz, 
\emph{Automorphisms of 
a free associative algebra of rank 2. I, II},
Trans. Amer. Math. Soc., I: 160 (1971),
393--401; II: 171  (1972), 309--315.
%
\bibitem[Dai19]{Dai19}
Daniel~Daigle,
\emph{Locally nilpotent sets of derivations}.
In: Kuroda, S., Onoda, N., 
Freudenburg, G. (eds). 
Polynomial Rings and 
Affine Algebraic Geometry, 41--71. 
PRAAG 2018. Springer Proceedings in 
Mathematics $\&$ Statistics, 
vol 319, 2019. Springer, Cham.
%
\bibitem[Dani89]{Dani89}
W\l{}odzimierz~Danielewski, 
\emph{On the cancellation problem 
and automorphism group of 
affine algebraic varieties}, 
Unpublished preprint, 1989.
%
\bibitem[Dan74]{Dan74}
Vladimir~ I.~Danilov, 
\emph{Non-simplicity of the group of 
unimodular automorphisms 
of an affine plane}, 
Mat. Zametki 15 (1974), 289--293.
%
\bibitem[Dem82]{Dem82}
Aleksander~S.~Demushkin, 
\emph{Combinatorial invariance of toric singularities}, 
Vestnik Moskov. Univ. Ser. I Mat. Mekh. 2 (1982), 
80--87.
%
\bibitem[D\'es06]{Des06}
Julie D\'eserti,
\emph{On the group of polynomial automorphisms 
of the affine plane},
J. Algebra 297:2 (2006), 584--599. 
%
\bibitem[Dub04]{Dub04}
Adrien~Dubouloz, 
\emph{Completions of normal affine surfaces 
with a trivial Makar-Limanov invariant}, 
Michigan Math. J. 52 (2004), 289--308. 
%
\bibitem[DMJP10]{DMJP10}
Adrien~Dubouloz,  Lucy~Moser-Jauslin, and
Pierre-Marie~Poloni,
\emph{Inequivalent embeddings of the
Koras--Russell cubic $3$-fold.}
Michigan Math. J. 59 (2010), 679--694.
%
\bibitem[DG15]{DG15}
Amartya~K.~Dutta and Neena~Gupta,
\emph{The epimorphism theorem and its generalizations},
J. Algebra Appl. 14:9 (2015), 1540010.
%
\bibitem[Edo18]{Edo18} \'Eric~Edo. 
\emph{Closed subgroups of the polynomial
automorphism group containing the
affine subgroup}. 
Transform.\ Groups 23 (2018), 71--74.
%
\bibitem[EL19]{EL19} 
\'Eric~Edo and Drew~Lewis, 
\emph{Co-tame polynomial automorphisms}, 
Internat.~J.~Algebra Comput.\  
29 (2019), 803--825. 
%
\bibitem[EP15]{EP15}
\'Eric~Edo and Pierre-Marie~Poloni,
\emph{On the closure of the tame automorphism group},
Int. Math. Res. Not. 19 (2015), 9736--9750. 
%
\bibitem[EGS23]{EGS23}
Polina Evdokimova, Sergey Gaifullin, 
and Anton Shafarevich,
\emph{Rigid trinomial varieties},
arXiv:2307.06672 (2023).
%
\bibitem[Fie94]{Fie94}
Karl-Heinz~Fieseler, 
\emph{On complex affine surfaces with 
$\CC^+$-action},
Comment. Math. Helv. 69:1  (1994),  5--27.
%
\bibitem[FKZ07]{FKZ07}
Hubert~Flenner, Shulim~Kaliman, 
and Mikhail Zaidenberg,
\emph{Completions of $\CC^*$-surfaces}. 
In: \emph{Affine algebraic geometry}, 149--200.
Osaka Univ. Press, Osaka 2007.
%
\bibitem[FKZ09]{FKZ09}
Hubert~Flenner, Shulim~Kaliman, 
and Mikhail~Zaidenberg,
\emph{On the Danilov-Gizatullin 
isomorphism theorem}, 
Enseign. Math. (2) 55:3-4 (2009), 275--283.
%
\bibitem[FKZ11]{FKZ11}
Hubert~Flenner, Shulim~Kaliman, 
and Mikhail~Zaidenberg,
\emph{Smooth affine surfaces with non-unique 
$\CC^*$-actions}, 
J. Algebraic Geom. 20 (2011), 329--398.
%
\bibitem[FKZ16]{FKZ16}
Hubert~Flenner, Shulim~Kaliman, 
and Mikhail~Zaidenberg,
\emph{The Gromov-Winkelmann 
theorem for flexible varieties}, 
J.\ Eur.\ Math.\ Soc.\ 18:11 
(2016), 2483--2510.
%
\bibitem[FZ05b]{FZ05b}
Hubert~Flenner and Mikhail~Zaidenberg, 
\emph{On the uniqueness of 
$\CC^*$-actions on affine surfaces}, 
In: Affine Algebraic Geometry,
Contemporary Mathematics Vol. 369, 97--111. 
Amer. Math. Soc., Providence, RI, 2005.
%
\bibitem[Fre17]{Fre17}
Gene~Freudenburg,
\emph{Algebraic Theory of Locally 
Nilpotent Derivations}, 
2nd edn. 
\emph{Encyclopaedia of Mathematical Sciences} 136. 
\emph{Invariant Theory and Algebraic 
Transformation Groups} 
VII. Berlin: Springer, 2017.
%
\bibitem[FM15]{FM15}
Gene~Freudenburg and Lucy~Moser-Jauslin,
\emph{A nonlinearizable action of 
$S_3$ on $\CC^4$},
Ann. Inst. Fourier, Grenoble, 52:1 (2002), 133--143.
%
\bibitem[Fur15]{Fur15}
Jean-Philippe~Furter, 
\emph{Polynomial composition rigidity and plane
polynomial automorphisms}, 
J.\ London Math.\ Soc.\ (2) 91 (2015), 180--202.
%
\bibitem[FuKr18]{FuKr18}
Jean-Philippe~Furter and Hanspeter~Kraft, 
\emph{On the geometry of automorphism
groups of affine varieties}, arXiv:1809.04175 (2018).
%
\bibitem[FuMa07]{FuMa07}
Jean-Philippe~Furter and  Stefan~Maubach, 
\emph{Locally finite polynomial endomorphisms},
 J. Pure Appl. Algebra 211 (2007), 445--458.
%
\bibitem[FuPo18]{FuPo18}
Jean-Philippe~Furter and Jean-Marie~Poloni, 
\emph{On the maximality of the 
triangular subgroup},
Ann. Inst. Fourier 68:1 (2018), 393--421.
%
\bibitem[Fur83]{Fur83}
Mikio~Furushima, 
\emph{Finite groups of polynomial 
automorphisms in $\CC^n$}, 
Tohoku Math. J. (2) 35 (1983), 415--424.
%
\bibitem[Giz71a]{Giz71a}
Marat~Kh.~Gizatullin,
\emph{Affine surfaces that are 
quasihomogeneous with respect
to an algebraic group}, 
Math. USSR Izv. 5 (1971), 754--769.
%
\bibitem[Giz71b]{Giz71b}
Marat~Kh.~Gizatullin,
\emph{Quasihomogeneous affine surfaces}, 
Math. USSR Izv. 5 (1971), 1057--1081. 
%
\bibitem[Giz18]{Giz18}
Marat~Gizatullin,
\emph{Two examples of affine homogeneous varieties},
Eur. J. Math. 4:3  (2018), 1035--1064.
%
\bibitem[GD75]{GD75}
Marat~Kh.~Gizatullin and Vladimir~I.~Danilov, 
\emph{Automorphisms of affine surfaces I}, 
Math. USSR Izv. 9:3 (1975), 493--534.
%
\bibitem[GD77]{GD77}
Marat~Kh.~Gizatullin and Vladimir~I.~Danilov, 
\emph{Automorphisms of affine surfaces II}, 
Math. USSR Izv. 11:1 (1977), 51--98.
%
\bibitem[Gub98]{Gub98}
Joseph~Gubeladze,
\emph{The isomorphism problem 
for commutative monoid rings},
J. Pure Appl. Algebra 129  (1998), 35--65.
%
\bibitem[GMM21]{GMM21}
Rajendra~V.~Gurjar, Kayo~Masuda, 
and Masayoshi~Miyanishi, 
\emph{Affine space fibrations},
Studies in Mathematics vol. 79, De Gruyter, 2021. 
%
\bibitem[Jel15]{Jel15}
Zbigniew Jelonek, 
\emph{On the group of automorphisms 
of a quasi-affine variety}, 
Math. Ann. 362:1-2 (2015), 569--578. 
%
\bibitem[Jun42]{Jun42}
Heinrich~W.~E.~Jung, 
\emph{\"{U}ber ganze birationale 
Transformationen der Ebene}, 
J. Reine Angew. Math. 184 (1942), 161--174.
%
\bibitem[Jur90]{Jur90}
Jerzy Jurkiewicz,
\emph{Linearizing some $\ZZ/2\ZZ$ actions on affine space},
Compos. Math. 76:1-2 (1990), 243--245.
%
\bibitem[Kal04]{Kal04}
Shulim~Kaliman,
\emph{Free $\CC^+$-actions on $\CC^3$ 
are translations}, 
Invent. Math. 156:1 (2004), 163--173. 
%
\bibitem[Kal09]{Kal09}
Shulim~Kaliman,
\emph{Actions of $\CC^*$ and $\CC^+$ 
on affine algebraic varieties},
in: \emph{Algebraic Geometry}, Seattle 2005.
D. Abramovich e.a.
(eds.), 629--654.
Proceedings of Symposia in Pure mathematics
vol. 80, Part 2. Amer. Math. Soc. 2009.
%
\bibitem[KML97]{KML97}
Shulim~Kaliman and Leonid~Makar-Limanov, 
\emph{On Russell-Koras contractible threefolds}, 
J. Algebraic Geom. 6 (1997), 247--268.
%
\bibitem[KKMLR97]{KKMLR97} 
Shulim~Kaliman, Mariusz~Koras, 
Leonid~Makar-Limanov, and Peter~Russell, 
\emph{$\CC^*$-actions on $\CC^3$ are linearizable}, 
Electron. Res. Announc. Amer. Math. Soc. 3 
(1997), 63--71.
%
\bibitem[KZ99]{KZ99}
Shulim~Kaliman and Mikhail~Zaidenberg, 
\emph{Affine modifications and 
affine hypersurfaces with a very transitive 
automorphism group}, 
Transform. Groups 4 (1999), 53--95.
%
\bibitem[KZ00]{KZ00}
Shulim~Kaliman and Mikhail~Zaidenberg, 
\emph{Miyanishi's characterization 
of the affine $3$-space 
does not hold in higher dimensions}, 
Ann. Inst. Fourier (Grenoble) 50 (2000), 
1649--1669.
%
\bibitem[Kam75]{Kam75}
Tatsuji~Kambayashi, 
\emph{On the absence of nontrivial separable 
forms of the affine plane}.
J. Algebra 35  (1975), 449--456.
%
\bibitem[Kam79]{Kam79}
Tatsuji~Kambayashi, 
\emph{Automorphism group of a polynomial ring 
and algebraic group action on an affine space}, 
J. Algebra, 60 (1979), 439--451.
%
\bibitem[KamRu82]{KamRu82}
Tatsuji~Kambayashi and Peter~Russell, 
\emph{On linearizing algebraic torus actions},
J. Pure Appl. Algebra 23:3 (1982), 243--250.
%
\bibitem[Kno91]{Kno91}
Friedrich~Knop, 
\emph{Nichtlinearisierbare Operationen halbeinfacher 
Gruppen auf affinen R\"{a}umen},
Invent. Math. 105:1  (1991), 217--220.
%
\bibitem[KoRu89]{KoRu89}
Mariusz~Koras and Peter~Russell, 
\emph{Codimension $2$ torus actions 
on affine $n$-space}. 
In: \emph{Group actions and invariant theory, 
Proc. Conf., Montreal/Can. 1988}.
Canad. Math. Soc. Conf. Proc. 10 (1989), 103--110.
%
\bibitem[KoRu97]{KoRu97}
Mariusz~Koras and Peter~Russell, 
\emph{Contractible threefolds and 
$\CC^*$-actions on $\CC^3$}, 
J. Algebraic Geom. 6 (1997), 671--695.
%
\bibitem[KPZ17]{KPZ17}
Sergei~Kovalenko, Alexander~Perepechko 
and Mikhail~Zaidenberg, 
\emph{On automorphism groups of affine surfaces}. 
In: \emph{Algebraic varieties and automorphism groups}, 
pp. 207--286.
Adv. Stud. Pure Math., vol. {\bf 75}, 
Math. Soc. Japan, Tokyo, 2017.
%
\bibitem[Kra89]{Kra89}
Hanspeter~Kraft, 
\emph{$G$-vector bundles and the Linearization Problem}. 
In: \emph{Group Actions and Invariant Theory}, 
CMS Conference Proceedings 
vol. 10, 111--123.
Amer. Math. Soc., Providence, RI, 1989.  
%
\bibitem[Kra96]{Kra96}
Hanspeter~Kraft, 
\emph{Challenging problems on affine $n$-space}, 
 S\'eminaire Bourbaki, Exp. No. 802, 5, Vol. 1994/95.
 Ast\'erisque 237 (1996),  295--317.
%
\bibitem[Kra17]{Kra17}
Hanspeter~Kraft, 
\emph{Automorphism groups of affine varieties 
and a characterization of affine $n$-space}, 
Trans. Moscow Math. Soc. 78 (2017), 171--186. 
%
\bibitem[KrPo85]{KrPo85}
Hanspeter~Kraft and Vladimir~L.~Popov, 
\emph{Semisimple group actions 
on the three-dimensional affine space are linear}, 
Comment. Math. Helv. 60 (1985), 466--479.
%
\bibitem[KrRe17]{KrRe17}
Hanspeter~Kraft and Andriy~Regeta, 
\emph{Automorphisms of 
the Lie algebra of vector fields on affine $n$-space},
J. Eur. Math. Soc. 19:5  (2017), 1577--1588.
%
\bibitem[KrRu14]{KrRu14}
Hanspeter~Kraft and Peter~Russell, 
\emph{Families of group actions, 
generic isotriviality, and linearization}, 
Transform. Groups 19:3 (2014), 779--792.
%
\bibitem[KrSc92]{KrSc92}
Hanspeter~Kraft and Gerald~W.~Schwarz, 
\emph{Reductive group actions 
with one-dimensional quotient}, 
Inst. Hautes Etudes Sci. Publ. Math. 76 (1992), 1--97.
%
\bibitem[KrSc95]{KrSc95}
Hanspeter~Kraft and Gerald~Schwarz, 
\emph{Finite automorphisms of affine $N$-space}. 
In: Essen, Arno van den (ed.), 
\emph{Automorphisms of affine spaces}. 
Proceedings of the international conference 
on invertible polynomial maps, 
held in Cura\c cao, Netherlands Antilles, July 4-8, 1994.
Dordrecht: Kluwer Academic Publishers. 55--66 (1995).
%
\bibitem[KrSt13]{KrSt13}
Hanspeter~Kraft and Immanuel~Stampfli, 
\emph{On automorphisms of the affine Cremona group},
Ann. Inst. Fourier 63:3  (2013), 1137--1148.
%
\bibitem[KrZa24]{KrZa24}
Hanspeter~Kraft and Mikhail~Zaidenberg,
\emph{Algebraically generated groups 
and their Lie algebras},
J. London Math. Soc. (2) 109 (2024), 1--39. 
\bibitem[Kul93]{Kul93}
Viktor~S.~Kulikov, 
\emph{Generalized and
local Jacobian problems}, 
Izvestiya
Mathematics 41:2 (1993), 351--365.
%
\bibitem[Kum02]{Kum02}
Shrawan~Kumar, 
\emph{Kac-Moody groups, 
their flag varieties and representation theory}, 
Progress in Mathematics, vol. 204, 
Birkh\"{a}user Boston Inc., Boston, MA, 2002. 
%
\bibitem[Lam01]{Lam01}
St\'ephane~Lamy, 
\emph{L’alternative de Tits pour $\Aut[\CC^2]$}, 
J. Algebra 239:2 (2001), 413--437. 
%
\bibitem[Lam24]{Lam24}
St\'ephane~Lamy,~
\emph{The Cremona group}.~Preliminary~version, 
July 27, 2024, 
available at: 
\href{https://www.math.univ-toulouse.fr/~slamy/blog/cremona.html}
{https://www.math.univ-toulouse.fr/~slamy/blog/cremona.html}
%
\bibitem[LP21]{LP21}
St\'ephane~Lamy and Piotr~Przytycki, 
\emph{Presqu'un immeuble pour 
le groupe des automorphismes mod\'er\'es}, 
Ann. H.~Lebesgue 4 (2021), 605--651.
%
\bibitem[LP22]{LP22}
St\'ephane~Lamy and Piotr~Przytycki, 
\emph{Tits alternative for the $3$-dimensional tame
automorphism group},
arXiv:2206.05611 (2022). 
%
\bibitem[LV13]{LV13}
St\'ephane~Lamy and St\'ephane~V\'en\'ereau,
\emph{The tame and the wild automorphisms 
of an affine quadric threefold},
J. Math. Soc. Japan 65:1 (2013), 299--320.
%
\bibitem[LRU23]{LRU23}
Alvaro~Liendo, Andriy~Regeta, and Christian~Urech, 
\emph{Characterization of affine surfaces 
with a torus action by their automorphism groups}, 
Ann. Sc. Norm. Super. Pisa, Cl. Sci. (5) 24:1  (2023), 
249--289.
%
\bibitem[LPU23]{LPU23}
Anne Lonjou, Piotr Przytycki, and Christian Urech,
\emph{Finitely generated subgroups 
of algebraic elements 
of plane Cremona groups are bounded},
arXiv:2307.01334 (2023).
%
\bibitem[Lun73]{Lun73}
Domingo~Luna,
\emph{Slices \'etales},
M\'emoires de la S. M. F. 33 (1973), 81--105.
%
\bibitem[ML70]{ML70}
Leonid~G.~Makar-Limanov, 
\emph{Automorphisms of a free algebra with 
two generators}, 
Funct. Anal. Appl., 4:3 (1970), 262--264.
%
\bibitem[Mas03]{Mas03}
Kayo~Masuda,
\emph{Nonlinearizable actions of dihedral groups 
on affine space},
Trans. Amer. Math. Soc. 356:9 (2003), 3545--3556.
%
\bibitem[MMP91]{MMP91}
Mikiyo~Masuda, Lucy~Moser-Jauslin,  and Ted~Petrie,  
\emph{Equivariant algebraic vector bundles 
over representations of
reductive groups: Applications}, 
Proc. Nat. Acad. Sci. 88 (1991), 9065--9066.
%
\bibitem[MMP96]{MMP96}
Mikiyo~Masuda, Lucy~Moser-Jauslin,  and Ted~Petrie,
\emph{The equivariant Serre problem for abelian groups}, 
Topology 35 (1996), 329--334. 
%
\bibitem[MP91]{MP91}
Mikiyo~Masuda and Ted~Petrie, 
\emph{Equivariant algebraic vector bundles 
over representations of reductive groups:
Theory}, 
Proc. Nat. Acad. Sci. 88 (1991), 9061--9064.
%
\bibitem[MP95]{MP95}
Mikiyo~Masuda and Ted~Petrie, 
\emph{Stably trivial equivariant algebraic 
vector bundles}, 
J. Amer. Math. Soc. 8 (1995), 687--714. 
%
\bibitem[Mat23]{Mat23}
Olivier~Mathieu,
\emph{Linearity and nonlinearity of groups of
polynomial automorphisms of the plane},
J.~Algebra 637 (2024), 47--89. 
%
\bibitem[MP09]{MP09}
Stefan~Maubach and Pierre-Marie~Poloni,
\emph{The Nagata automorphism is 
shifted linearizable},
J. Algebra 321 (2009), 879--889.
%
 \bibitem[Miy24]{Miy24}
Masayoshi~Miyanishi, 
\emph{Affine algebraic geometry: 
geometry of polynomial rings}.
Series on University Mathematics, vol. 11. 
World Scientific, New Jersey, 2024. 
%
\bibitem[MJ94]{MJ94}
Lucy~Moser-Jauslin,
\emph{Algebraic equivariant vector bundles
and the linearity problem}, 
in: 
Proc. Symp. in Pure Mathem.
vol. 56, Part I, 1994, 355--364.
%
\bibitem[Nag59]{Nag59}
Hirosi~Nagao, 
On $\GL(2, K[x])$, 
J. Inst. Polytech. Osaka City Univ. Ser. A 10 (1959),
117--121.
%
\bibitem[Nag72]{Nag72}
Masayoshi~Nagata, 
\emph{On automorphism group of $k [x, y]$},
Kinokuniya Book-Store Co., Ltd. Tokyo, 
Japan, 1972.
%
\bibitem[Pan84-86]{Pan84-86}
Dmitri~I.~Panyushev,
\emph{Semisimple automorphism groups
of four-dimensional affine space},
Math. USSR Izv. 23 (1984), 171--183.
Correction: \emph{ibid}. 27:3 (1986), 607--608.
%
\bibitem[Per24]{Per24}
Alexander~Perepechko, 
\emph{Structure of connected nested 
automorphism groups}, 
arXiv:2312.08359v2 (2024).
%
\bibitem[PR24]{PR24}
Alexander~Perepechko and Andriy~Regeta,
\emph{Automorphism groups of affine varieties 
without non-algebraic elements},  
Proc. Amer.Math. Soc.
152 (2024), 2377--2383 (to appear). 
%
\bibitem[PR23]{PR23}
Alexander~Perepechko and Andriy~Regeta,
\emph{When is the automorphism 
group of an affine variety nested?},
Transform. Groups 28:1 (2023), 401--412. 
%
\bibitem[PZ22]{PZ22}
Alexander~Perepechko and 
Mikhail~Zaidenberg, 
\emph{Automorphism groups of 
rigid affine surfaces: 
the identity component}, 
arXiv:2208.09738 (2022).
%
\bibitem[PeRa86]{PeRa86}
Ted~Petrie and John~D.~Randall, 
\emph{Finite-order algebraic automorphisms 
of affine varieties}, 
Comment. Math. Helv. 61:2 (1986), 203--221. 
%
\bibitem[Pop73]{Pop73}
Vladimir~L.~Popov, 
\emph{Classification of affine algebraic surfaces 
that are quasihomogeneous with respect to 
an algebraic group}, 
Math. USSR Izv. 7 (1973),
1039--1055 (1975).
%
\bibitem[Pop87]{Pop87}
Vladimir~L.~Popov, 
\emph{On actions of $\Ga$ on $\AA^n$}, 
in: \emph{Algebraic groups. {U}trecht 1986}, 237--242.
Lecture Notes in Math., vol. 1271, Springer, Berlin,
1987. 
%
\bibitem[Pop14]{Pop14}
Vladimir~L.~Popov,
\emph{On infinite dimensional algebraic transformation groups}, 
Transform. Groups 19 (2014), 549--568.
%
\bibitem[Pop17]{Pop17}
Vladimir~L.~Popov,
\emph{Bass' triangulability problem},
in: \emph{Algebraic varieties and automorphism groups}, 
425--441.
Adv. Stud. Pure Math. vol. 75, 
Math. Soc. Japan, Tokyo, 2017.
%
\bibitem[Smi89]{Smi89}
Martha~K.~Smith,
\emph{Stably tame automorphisms},
 J. Pure Appl. Algebra {\bf 58} (1989), 209--212.
%
\bibitem[Ram64]{Ram64} Chakravarthi~P.~Ramanujam, 
\emph{A note on automorphism 
groups of algebraic varieties}, 
Math. Ann. 156 (1964), 25--33.
%
\bibitem[RvS21]{RvS21}
Andriy~Regeta and Immanuel~van~Santen,
\emph{Maximal commutative unipotent 
subgroups and a characterization 
of affine spherical varieties}, 
J. Eur. Math. Soc. (to appear), 
arXiv:2112.04784 (2021). 
%
\bibitem[Ren68]{Ren68}
Rudolf~Rentschler, 
\emph{Op\'erations du groupe additif sur le plan affine}, 
C. R. Acad. Sci. Paris Ser. A 267 (1968), 384--387.
%
\bibitem[RCS04]{RCS04}
Vitalij~A.~Roman'kov, Igor~V.~Chirkov, 
and Mikhail~A.~Shevelin, 
\emph{Nonlinearity of the automorphism groups 
of some free algebras}, 
Siberian Math. J. 45:5 (2004), 974--977.
%
\bibitem[Rud86]{Rud86}
Alexei~N.~Rudakov, 
\emph{Subalgebras and automorphisms 
of Lie algebras of Cartan type}, 
Func. Anal. Appl. 20:1 (1986), 72--73.
%
\bibitem[Sch89]{Sch89}
Gerald~W.~Schwarz, 
\emph{Exotic algebraic group actions}, 
C. R. Acad. Sci. Paris Ser. Math. 309 (1989), 89--94.
%
\bibitem[Ser03]{Ser03}
Jean-Pierre~Serre,
\emph{Trees},
Springer Monographs in Mathematics. 
Springer-Verlag, Berlin, 2003.
%
\bibitem[Sha81]{Sha81}
Igor~R.~Shafarevich, 
\emph{On some infinite-dimensional groups}. II, 
Izv. Akad. Nauk SSSR Ser. Mat. 45:1 (1981), 214--226.
%
\bibitem[Sha04]{Sha04}
Igor~R.~Shafarevich, 
\emph{On the Group $\GL(2,K[t])$},
Proc. Steklov Inst. Math. 246  (2004), 308--314. 	
%
\bibitem[SU04a]{SU04a} 
Ivan~P.~Shestakov and Ualbai~U.~Umirbaev, 
\emph{Poisson brackets and two generated 
subalgebras of rings of polynomials},
 J. Amer. Math. Soc. 17 (2004), 181--196. 
 %
 \bibitem[SU04b]{SU04b} 
Ivan~P.~Shestakov and Ualbai~U.~Umirbaev, 
\emph{Tame and wild automorphisms of rings 
of polynomials in three variables}, 
 J. Amer. Math. Soc. 17 (2004), 197--227. 
%
\bibitem[SZ24]{SZ24}
Ivan~Shestakov and Zerui Zhang,
\emph{An Anick type wild automorphism 
of free Poisson algebras},
arXiv:2407.04919 (2024).
%
\bibitem[SY98]{SY98}
Vladimir~Shpilrain and Jie-Tai~Yu,
 \emph{On generators of polynomial algebras 
 in two commuting or non-commuting variables},
J. Pure Appl. Algebra 132 (1998), 309--315.
%
\bibitem[Sku21]{Sku21}
Alexander~A.~Skutin,
\emph{Maximal Lie subalgebras 
among locally nilpotent derivations},
Sb. Math. 212:2 (2021), 265--271.
%
\bibitem[Smi89]{Smi89}
Martha~K.~Smith, 
\emph{Stably tame automorphisms}, 
J. Pure Appl. Algebra 58:2 (1989), 209--212. 
%
\bibitem[Sno89]{Sno89}
Denis~Snow, 
\emph{Unipotent actions on affine space}, 
in: \emph{Topological Methods in 
Algebraic Transformation Groups 
(New Brunswick, NJ, 1988)}, 
165--176.
Progr. Math., vol. 80, Birkhauser Boston, 
Boston, MA, 1989.
%
\bibitem[Sos07]{Sos07}
Yury~V.~Sosnovskii, 
\emph{The hypercentral structure of the group of 
unitriangular automorphisms of a polynomial algebra}, 
Siberian Math. J. 48:3 (2007), 555–-558.
%
\bibitem[Spo07]{Spo07}
Stanis\l{}aw~Spodzieja,
\emph{On the Nagata automorphism},
Univ. Iagel. Acta Math. XIV (2007), 131--136.
%
\bibitem[Sta13]{Sta13}
Immanuel~Stampfli, 
\emph{A note on automorphisms of the affine Cremona group},
Math. Res. Lett. 20:6  (2013), 1177--1181.
\bibitem[Umi07]{Umi07}
Ualbai~U.~Umirbaev,  
\emph{The Anick automorphism of free associative algebras},
J. Reine Angew. Math. 605 (2007), 165--178.
%
\bibitem[Ure21]{Ure21}
Christian~Urech, 
\emph{Subgroups of elliptic elements 
of the Cremona group}, 
J. Reine Angew. Math. 770 (2021), 27--57. 
%
\bibitem[UZ21]{UZ21}
Christian~Urech and Susanna~Zimmermann,
\emph{Continuous automorphisms of Cremona groups},
Int. J. Math. 32:4 (2021), Article ID 2150019, 17 p. 
%
\bibitem[vdE00]{vdE00}
Arno~van~den~Essen, 
\emph{Polynomial automorphisms 
and the Jacobian conjecture}, 
Progress in Mathematics, vol. 190, 
Birkh\"{a}user, Basel, Boston, Berlin, 2000.
%
\bibitem[vdEKC21]{vdEKC21}
Arno~van~den~Essen, Shigeru~Kuroda, 
and Anthony~J.~Crachiola,
\emph{Polynomial automorphisms 
and the Jacobian conjecture. 
New results from the beginning
of the 21st Century}, in: \emph{Frontiers in Mathematics}. 
Birkh\"{a}user, Cham, Switzerland, 2021.
%
\bibitem[vdEMV07]{vdEMV07}
Arno~van~den~Essen, Stefan~Maubach,
and Stephane~V\'en\'ereau, 
\emph{The special automorphism group of 
$R[t]/(t^m)[x_1, \ldots, x_n]$ and
coordinates of a subring of $R[t][x_1, \ldots , x_n]$},
J. Pure Appl. Alg. 210 (2007), 141--146.
%
\bibitem[vdK53]{vdK53}
Wouter~van~der~Kulk, 
\emph{On polynomial rings in two variables}, 
Nieuw Arch. Wisk. (3) 1 (1953), 33--41.
%
\bibitem[Win90]{Win90}
J\"{o}rg~Winkelmann, 
\emph{On free holomorphic $\CC_+$-actions 
on $\CC^n$ and homogeneous Stein manifolds}, 
Math. Ann. 286:1-3 (1990),  593--612.
%
\bibitem[Wri15]{Wri15}
David~Wright,
\emph{The generalized amalgamated product structure 
of the tame automorphism group in dimension three},
Transform. Groups 20:1 (2015), 291--304.
%
\bibitem[Wri17]{Wri17}
David~Wright,
\emph{Amalgamation and automorphism groups},
in: \emph{Algebraic varieties and automorphism groups}, 
465--474.
Adv. Stud. Pure Math. vol. 75, 
Math. Soc. Japan, Tokyo, 2017.
%
\bibitem[Zyg11]{Zyg11}
Jakub~Zygadło,
\emph{Remarks on a normal subgroup of $GA_n$}, 
Comm. Algebra 39:6 (2011), 1992--1996. 
%
\end{thebibliography}
\end{document}